\documentclass[12pt]{article}
\usepackage{amsmath}
\usepackage{amssymb}
\usepackage{amsthm}

\newtheorem{theorem}{Theorem}
\newtheorem{lemma}{Lemma}
\newtheorem{note}{Note}
\newtheorem{corollary}{Corollary}

\newtheorem{condition}{Condition}
\newtheorem{assumption}{Assumption}

\newtheorem{definition}{Definition}
\newcommand{\diag}{\mathop{\rm diag}\nolimits}

\renewcommand{\det}{\mathop{\rm det}\nolimits}
\newcommand{\tr}{\mathop{\rm tr}\nolimits}
\newcommand{\adj}{\mathop{\rm adj}\nolimits}
\def\Proof{{\bf Proof}}

\font\Bbb=msbm10 scaled 1200
\def\R{\hbox{\Bbb R}}
\def\S{\hbox{\Bbb S}}

\def\C{\hbox{\Bbb C}}
\def\H{\hbox{\Bbb H}}
\def\I{\hbox{\Bbb I}}
\def\P{\hbox{\Bbb P}}
\def\beq{\begin{equation}} \def\eeq{\end{equation}}
\font\gothic=eufm10
\def\goth#1{\hbox{\gothic #1}}

\begin{document}

\begin{center}
          FIELD THEORY FOR MULTIPLE INTEGRALS \\
                       M.I.Zelikin       \\

\end{center}

\bigskip

 \begin{abstract}

 New constructions in the theory of fields for multiple integrals
 are designed. Generalizations of the Legendre - Weyl -
 Caratheodory transforms and corresponding invariant integrals are
 introduced and explored. Connection and curvature of bundles induced
by a field of extremals are calculated.

 \end{abstract}

\section{Introduction}

Let $\goth N$ be a domain on a smooth $n$-dimensional Riemannian
manifold, and let $\rho : \xi \to \goth N$ be a $\nu$-dimensional
vector bundle over the base $\goth N$. A fiber of the bundle over
a point $t \in \goth N$, i.e. the full inverse image of the point
$t$ under the map $\rho$, is an $\nu$-dimensional linear subspace.
Local coordinates on $\goth N$ are denoted by $t = (t^1,...t^n)$;
local coordinates on fibers are denoted by $x = (x^1,...x^{\nu})$.
Usual convention on summation over two-fold occurring indices is
used. Latin indices correspond to coordinates on the base and vary
from $1$ to $n$, Greek ones correspond to coordinates on fibers
and vary from $1$ to $\nu$. Multi-indices will be denoted by
capital Latin and Greek indices respectively. Collection of all
indices from $1$ to $n$ is denoted by $I$. The ordered exterior
product of differential entering into the multi-index $K$ is
denoted by $dt^K := dt^{i_1}\wedge dt^{i_2} \wedge ... dt^{i_k}$
(the letter $K$ shows that $|K|=k$, where the absolute value of
multi-index means the number of its indices). The symbol wedge for
exterior product (the symbol $\wedge$) will be omitted for the
sake of brevity; it always will be implied while dealing with
product of differentials.

Consider the functional whose part related to the chart
 $V \subset \goth N$, is

\beq {\cal F} = \int_V f \left ( t,x,\frac {Dx}{Dt} \right ) dt^I.
 \label{1}
 \eeq

Subsequent calculations will be produced in coordinates of the
chart $V$.

Let us denote by $J_1(\xi)$ the bundle of 1-jets over $\xi$ and
let
 \beq
 q^{\alpha}_i := \frac {\partial x^{\alpha}}{\partial
t^i}= g^{\alpha}_i(t,x)
 \label{2}
 \eeq

\noindent be a section of $J_1(\xi)$. Such section can be
considered as a slope field $\cal G$, i,e. as a distribution of
$n$-dimensional planes in the space $\xi$. We shall say that a
manifold

 $M=\{x=\hat x(t)\}$  embedded in the slope field
$\cal G$, if $\frac {d\hat x^\alpha}{dt^i} = g^{\alpha}_i(t,x)$.

We distinguish three standpoints concerning the arguments in our
functions:

1. If $t^i,x^\alpha,\frac {\partial x^{\alpha}}{\partial t^i}$ are
taken as independent variables, as for instance in the function
$L$, then the derivatives with respect to these variables are
marked by attaching the respective variable as an index.

2. By using a given slope field $\cal G$, the arguments $\frac
{d\hat x^\alpha}{dt^i}$ are replaced by functions
$g^{\alpha}_i(t,x)$ that depends on the $t^i$ and $x^\alpha$ only.
The partial derivatives with respect to the arguments $t^i$ and
$x^\alpha$ are then denoted by $\frac {\partial }{\partial t^i}$
and $\frac {\partial }{\partial x^{\alpha}}$.

3. The substitution $x(t)$ and its derivatives referring to a
given surface $V$ changes functions which appeared in the second
(or the first) standpoint, into functions of the $t$ alone. Their
derivation with respect to $t^i$ is denoted by $\frac {d}{dt^i}$.

The vanishing of the first variation is expressed by Euler's
equations

$$
\frac {df_{\frac {\partial x^\alpha}{\partial t^i}}}
{dt^i} -
f_{x^\alpha} = 0.
$$

Functions $x(\cdot)$ which satisfy this equation are called extremals.
The Jacobi matrix of variables $S^K$ relative to the arguments $t^J$
is denoted by $\frac {{\cal D}[S^K]}{{\cal D}[t^J]}$. The
calligraphical ${\cal D}$ means the derivative with respect to the
explicitly entering argument, while the direct $D$ (say $\frac
{D[S^K]}{D[t^J]}$) means the full derivative taking into account the
dependance $x(t)$. The determinant of this matrix is denoted by
 $|\frac {D[S^K]}{D[t^J]}|$. The identity matrix is denoted by
 $\I$ (its dimension is implicitly defined by the formula in question).

Nonnegativity of the second variation is the natural necessary
condition of minimum~\cite{Fo}. The investigation of conditions of
the nonnegativity begins with works of A.Clebsch~\cite{Cl} who
explored the Dirichlet functional

$$
\delta^2 {\cal F} =\int_V \left [ \frac {\partial^2 \hat f}
{\partial \left (\frac {\partial x^\alpha}{\partial t^i}\right)
\partial \left (\frac {\partial x^\beta}{\partial t^j}\right )}
\frac {\partial h^\alpha}{\partial t^i} \frac {\partial
h^\beta}{\partial t^j} + 2\frac {\partial^2 \hat f} {\partial \left
(\frac {\partial x^\alpha}{\partial t^i}\right) \partial x_\beta}
\frac {\partial h^\alpha}{\partial t^i} h^\beta + \frac
{\partial^2 \hat f}{\partial x^\alpha \partial x^\beta}h^\alpha
h^\beta \right ] dt^I.
$$

The idea of Clebsch was to reduce this functional to the integral
from its main part, i.e. to the rearrange quadratic terms relative
to first derivatives. The reduction was realized by adding under
the integral sign a closed differential form having the type of
divergency. It seems that Clebsch presumes that for multiple
integral a direct analog of Legendre condition is valid: if the
second variation is nonnegative then the quadratic terms relative
to first derivatives that is defined on the space of $(n \times
\nu)$-matrices
$$
q^\alpha_i = \frac {\partial x^\alpha}{\partial t^i}
$$
must be nonnegative for all values of $t$.

But J.Hadamard~\cite{H} almost through half a year after the work
of Clebsch shows that it is not correct. This quadratic form is
nonnegative not for all matrices. The correct necessary condition
(which is called Hadamard-Legendre condition) is the following.

\begin{theorem}[Hadamard]

Let the functional

\beq
\delta^2 {\cal F} =\int_V \left [a^{ij}_{\alpha \beta}(t)
\frac {\partial x^\alpha}{\partial t^i} \frac {\partial
x^\beta}{\partial t^j} + 2c^i_{\alpha \beta} (t) \frac {\partial
x^\alpha}{\partial t^i}x^\beta + b_{\alpha \beta}(t)x^\alpha
x^\beta \right ]dt^I
 \label{3}
 \eeq
be nonnegative for functions $x(\cdot)$ that meet boundary
conditions

$$
x|_{\partial V} = 0.
$$

Then for all values of $t \in V$ the quadratic form
$a^{ij}_{\alpha \beta}(t)q^\alpha _i q^\beta _j$ takes nonnegative
values on $(n \times \nu)$-matrices having the form $q^\alpha _i =
\xi^\alpha \eta _i$ (that is for matrices of rank 1)
\end{theorem}

The assertion of the theorem can be reformulated as follows: the
biquadratic form $a^{ij}_{\alpha \beta}(t)\xi^\alpha \xi^\beta
\eta _i\eta _j$ is nonnegative for all values of $t \in V$ and
$\xi \in \R^\nu, \eta \in (\R^n)^*$.

Let us remark that a simple sufficient condition for optimality of
small pieces of extremals is the condition of convexity of $f$
relative to variables $Dx/Dt.$ But this assumption is much more
strong than the necessary condition of Hadamard-Legendre.

The work of Hadamard stimulates whole series of works that aim for
decreasing the gap between necessary and sufficient conditions for
optimality. One of the bright work of this cycle was that of
Terpstra~\cite{T} where abstract algebraic questions suggested by
this themes were considered.

Namely, let

\beq
 a^{ij}_{\alpha \beta}(t)q^\alpha _i q^\beta _j
  \label{4}
 \eeq
be a quadratic form on $(n \times \nu)$-matrices $\| q^\alpha
_i\|$. Consider the cone of rank 1 matrices, that is ones having
the form $q^\alpha _i = \xi^\alpha \otimes \eta_i$. On this cone
the form~(\ref{4}) turns into a biquadratic form

\beq
a^{ij}_{\alpha \beta}(t)\xi^\alpha \xi^\beta \eta _i\eta_j,
\label{5}
 \eeq
which is defined on pairs of vectors $\xi \in \R^\nu, \, \eta \in
(\R^n)^*$. Assume that the form~(\ref{4}) is nonnegative on
matrices of rank 1. Terpstra set up the question: is it possible
to turn~(\ref{4}) into a positive form on the space of all
matrices by adding terms

$$
r^{ij}_{\alpha \beta}(q^\alpha_i q^\beta_j - q^\beta_i
q^\alpha_j),
$$
where tensor $r^{ij}_{\alpha \beta}$ is skew-symmetric relative to
both $i,j$, and $\alpha ,\beta$? Similar additions to the main
terms arise by adding a closed differential form to the integrand.
These additional terms give zero on rank 1 matrices, and the
biquadratic form~(\ref{5}) do not changes. Terpstra shows that it
is possible under the following condition

\begin{condition}
The form~(\ref{5}) admits decomposition into a sum of squares of
bilinear forms from which as a minimum $n\nu$ forms are linearly
independent.
\end{condition}

The question on a decomposition of even forms into a sum of
squares was investigated by Hilbert. The final result is due to
Van-der-Waerden who proves that a decomposition is realizable if
$min (\nu , n) \le 2$. If this minimum is more or equal to $3$
then there exist indecomposable positively definite forms. The
first explicit example of such forms was constructed by
Terpstra~\cite{T}. More simple example was suggested later by
D.Serre~\cite{S}, \cite{SE}.

This subject appears closely related with problems of existence of
minima. The leading part in the proofs of existence plays the
condition of lower semicontinuity of a functional in question.

C.Morrey~\cite{Mo} and later J.Ball~\cite{Ba} show that necessary
and sufficient condition for lower semicontinuity of second
variation is the Hadamard-Legendre condition. And they explored an
interesting generalization of the notion of convexity for multiple
integral that is called polyconvexity. To define it, consider an
integrand $f(t,x,\dot x)$ as a function of matrix
$$
q = \left \| \frac {\partial x^\alpha}{\partial t^i}\right \|
$$
where $t$ and $x$ are fixed. One corresponds to each $(n \times
\nu)$-matrix the set of elements of its exterior powers, that is all
 $(l \times l)$-minors of the matrix $1 \le l \le \min (n,\nu)$.
One obtains a point of $r$-dimensional space denoted by $\tau
(q)$. It is easy to calculate that $r = {n+\nu\choose n} - 1$. The
mapping $\tau$ transfers $\R^{n\nu}$ into an algebraic subset
${\cal K}$ of the space $\R^r$ which is defined by Pl\"ucker
relations on minors of the matrix $q$. Thus the function $f$
appears to be defined on the Pl\"ucker cone ${\cal K}$ of the
space $\R^r$. The function $f$ is called polyconvex if it admits a
convex extension on the whole space $\R^r$. The polyconvexity of
$f$ is a sufficient (but not a necessary) condition of lower
semicontinuity of the integral (see \cite{Mo}, \cite{Gi},
\cite{Sil}, \cite{Ba1}). It is interesting that functions
constructed by Terpstra and Serre gives examples of lower
semicontinuous but not polyconvex functionals. The gap between the
necessary and the sufficient conditions for optimality was
essentially shorten by Van-Hove~\cite{V} who proved that natural
strengthening of the Hadamard-Legendre condition

\beq
 \frac {\partial^2 \hat f} {\partial
\left (\frac {\partial x^\alpha}{\partial t^i}\right) \partial
\left (\frac {\partial x^\beta}{\partial t^j}\right )} \xi^\alpha
\xi^\beta \eta^i \eta^j \ge \varepsilon |\xi|^2|\eta|^2
 \label{6}
 \eeq
gives the locally sufficient condition for $C^1$-minimum. The
expression "locally sufficient" means that the domain of the
integration is sufficiently small.

The idea of the Van-Hove's proof is the following. First, we make
coefficients to be frozen, i.e. we fixe arguments
$(t=t_0,x=x_0)$ in coefficients of the quadratic form of the
integrand~(\ref{3}). It does not affects on the estimations since
it can be taken as a domain of integration an arbitrarily small
neighborhood of the point $(t=t_0,x=x_0)$. Second, we apply the
Fourier transform to the functional obtained and use Parseval
equality.

This construction disclose the internal reason why the
Hadamard-Legendre condition includes only rank $1$ matrices.
Namely, Fourier transform transfers the operation of
differentiation into the operation of multiplication by the
corresponding independent variable. So, if Fourier image of a
function $x^\alpha (t)$ is $\xi^\alpha (\eta)$ then the image of
its derivative $\frac {\partial x^\alpha}{\partial t^i}$ will be
$\xi^\alpha \eta_i$. As a result, the  biquadratic form which
stand at the left side of the formula~(\ref{6}) appears under the
integral sign of the Parseval equality, and the inequality~(\ref{6})
guarantee the positive definiteness of the functional in question.

To prove optimality conditions on "large" parts of the manifold,
one needs the theory of index of the functional, and the
generalization of the notion of the conjugate point. It was the
subject of big series of works (see, for example,
Dennemeyer~\cite{Den}, Simons~\cite{Si}, Smale~\cite{Sm},
Swanson~\cite{Sw}, \"Uhlenbeck~\cite{U}). We do not concern this
theme here.

The positivity of the second variation is insufficient to prove
the sufficient conditions of the strong minimum. One needs the
approach connected with the field theory and invariant Hilbert
integral~\cite{G}, \cite{F}. The first variant of this theory was
suggested by C.Caratheodory~\cite{C}. Another (relatively more
simple) variant with more strong demands on functions was
constructed by H.Weyl~\cite{W}. We will shortly describe both
approaches in a suitable form for subsequent presentation.

In variational calculus for multiple integrals one needs as many
principal functions of Hamilton (action-functions) as there are
independent variables. Let we have $n$ action-functions:
$S^i(t,x), \; (i = 1,...n)$. The Weyl construction is based on
invariant integral using divergence of a vector $S^i$:

$$
\Omega = \sum_{i=1}^n dt^1...dt^{i-1}dS^idt^{i+1}... dt^n = \frac
{dS^i}{dt^i} dt^I.
$$

Let us recall that under the expression $\frac {dS^i}{dt^i}$ we
mean the full derivative of the function $S^i$ taking into account
the dependence $x=x(t)$. For functions of matrices $q^\alpha_i =
\partial x^\alpha /\partial t^i$  Weyl suggests to use the direct
analog of Legendre transform in which the part of the scalar
product plays the trace of the product of matrices.

$$
       f(t,x,q) \mapsto f^*(t,x,p^*) = -f + \tr (p^*q), \quad
    (p^*)^i_\alpha =\frac {\partial f}{\partial q^\alpha _i}.
$$

Here the value $q$ as a function of $p^*$ has to give global
maximum to the Weierstrass function $-f + \tr (p^*q)$. So one
needs the condition of convexity of the function $f$ on the space
of matrices $q$ with fixed $t$ and $x$. The condition for the
minimum of the Weierstrass function is an equality

$$
                 \frac {\partial S^i}{\partial x^\alpha} =
\frac {\partial f}{\partial \left ( \frac {\partial x^\alpha}
{\partial t^i} \right )}.
$$

The construction of Caratheodory is based on an invariant integral
of the determinant from derivatives of functions $S^i$.

$$
\Omega = \det \left \| \frac {\partial S^i}{\partial t^j}
 \right \| dt^I.
$$

Caratheodory considered the Legendre transform as an algebraic
mapping of the space of quadruples. Let $(f,\phi ,q_i^\alpha
,p^i_\alpha ) \; 1\le i \le n, \; 1\le \alpha \le \nu$ be a set,
consisting of $2(n\nu + 1)$ elements, where $f$ is the initial
function, $\phi$ is the dual function, $q_i^\alpha$ are the
initial independent variables that correspond to
"velocities"($\partial x^\alpha/\partial t^i $), and $p^i_\alpha$
are the dual variables that correspond to moments $(\partial
f/\partial (\frac {\partial x^\alpha}{\partial t^i}))$. It is
supposed that all these variables are bound together by the relation

\beq
 f +\phi = \tr (pq).
 \label{7}
\eeq

We shall speak that a transformation is tangential if it transfers
the canonical differential form $df-pdq$ on the bundle of 1-jets
into proportional to it canonical form $df^* - q^*dp^*$.

New variables of the standard Legendre transform are obtained from
the old ones by simple permutation

$$
f^* = \phi, \; \phi^* = f, \; q^* = p, \; p^* = q.
$$

It is evident that this transformation is tangential, birational,
involutory, and it preserves the relation~(\ref{7}).

Caratheodory has constructed a new transformation in the space of
quadruples for the theory of invariant integrals having the type
of determinant. The original text of Caratheodory is difficult,
because of lack of motivations. Big series of many cumbersome
formulas with numerous indices, which at the end bring to desirable
results, give the impression of magic. So, we will give here the proof of
the main theorem of Caratheodory in convenient for us matrix-form.

Let us introduce an auxiliary square matrix $A$

\beq
 A=||f\I -  pq||.
\label{8}
 \eeq

Introduce the notation $\frac {f^{n-2}}{\det A} = \gamma$. The
transpose of a matrix $B$ will be denoted by $B^t$.

Define the mapping ${\cal Z}$ of a quadruple $(f,\phi,p,q)$ by the
formulae

\beq
f^* = \gamma f; \; \phi^* = \gamma \phi ; \ ; (p^*)^t =
\gamma qA; \; (q^*)^t=A^{-1}p.
 \label{9}
  \eeq

It follows
$$
p^*q^*= ApqA^{-1}\gamma .
$$
Using permutability of matrices $A$ and $pq$, one has
$$
 \frac {p^*q^*}{f^*}=\frac {pq}{f}
$$
Hence, taking into account the condition $f+\phi =\tr(pq)$, one
obtains
$$
f^* + \phi^* = \tr (p^*q^*).
$$

\smallskip

\begin{theorem}[Caratheodory].

The transformation ${\cal Z}$  is tangential, birational,
involutory, and it preserves the relation~(\ref{7}).
\end{theorem}

\Proof.

Rewrite the formula~(\ref{8}) as $pq = f\I - A$ and express
$(q^*)^t q$ using~(\ref{9}). One has  $(q^*)^t q = A^{-1}pq =
A^{-1}(f\I - A) = A^{-1}f - \I$. Hence, $\I + (q^*)^t q =
fA^{-1}$. It follows

\beq
 \det (\I + (q^*)^t q) = \frac {f^n}{\det A} = ff^*.
 \label{10}
 \eeq

Calculate the differential of~(\ref{10}) using the following
formula for differential of determinant

$$
d(\det g) = (\det g)\tr(g^{-1}dg).
$$

We have

\beq
\begin{array}{c}
fdf^* + f^*df = ff^* \left ( \tr(\frac {A}{f}[(dq^*)^tq +
(q^*)^tdq]
    \right ) = \\
f^*\left ( \tr (\frac {qA(dq^*)^t}{\gamma}) +
    \tr (A(q^*)^tdq) \right ) = \\

ff^* (\tr \left (\frac {p^*dq^*}{\gamma f}
    \right ) + \tr \left (\frac {pdq}{f} \right )) = \\

    f\tr (p^*dq^*) + f^*\tr (pdq).
\end{array}
\label{11}
 \eeq

By combining the first and the last members in the
chain~(\ref{11}), one has

\beq
 f(df^* - \tr (p^*dq^*) + f^*(df - \tr (pdq)) = 0.
  \label{12}
   \eeq

\smallskip

Now let us prove that the transformation ${\cal Z}$ is birational
and involutary. Indeed. Repeat the transformation ${\cal Z}\circ
{\cal Z}$. Introduce the auxiliary matrix $A^* =f^*\I - p^*q^*$.
Using~(\ref{8}) one obtains

$$
A^* =\gamma f\I - (\gamma qA)^t(A^{-1}p)^t = \gamma ||f\I -
(A^{-1}pqA)^t||.
$$

The matrix $pq$ commute with the matrix $A$, hence

\beq
 A^* = \gamma ||f\I - (pq)^t|| = \gamma A^t, \quad
 \gamma^* = \frac {(f^*)^{n-2}}{\det A^*} =
 \frac {\gamma^{n-2}f^{n-2}}{\gamma^n\det  A} = 1/\gamma .
\label{13}
 \eeq

Using~(\ref{9}), (\ref{13}) successively obtain

\beq
\begin{array}{c}
 f^{**} =\gamma^*(f^*) = f; \; \phi^{**} = \gamma^*
\phi^* = \phi; \\

 p^{**}=\gamma^*(q^*A^*)^t = \gamma^*\gamma AA^{-1}p = p; \\
  q^{**}=((A^*)^{-1}p^*)^t = \gamma (A^tq^t)^t\frac {1}{\gamma}A^{-1} = q.
\end{array}
 \label{14}
\eeq

$\Box$

\bigskip

Later on Th.De Donde~\cite{D}, J.T.Lepage~\cite{L},
H.B\"orner~\cite{B} a.o. connect approaches of H.Weyl and
C.Caratheodory with the technique of differential forms. They
tried to find the most general differential forms giving the field
theory.

We think that one has to find only invariant differential forms
having direct geometrical meaning (such as forms of Weyl and
Caratheodory). Action-functions define simultaneous
parametrization of all the manifolds which define "the flow" of
solutions generated a field of extremals. This is the reason why
differential forms that characterize "the flow" in the bundle
$\xi$ have to be invariant relative to choice of coordinates for a
mapping $\xi \to S$. Besides, it is not imperative that unknown
manifolds admit one-to-one projection onto the space of variables
$(t^1,...t^n)$. This fact counts in favour of invariance.

The main goal of this work is to find and to explore new
constructions of invariant integrals and the corresponding
tangential transformations. It will be designed a series of
transformations using tangent planes of $f(Dx/Dt)$ in spaces of
different exterior powers of matrices of jets. Note that
Caratheodory uses only one (the highest) exterior power --- the
determinant.

To construct a field theory one has to add to an integrand a
closed differential form, that does not changes the value of the
integral but turns the integrand into nonnegative function (analog
of the Weierstrass function). As in the theory of characteristic
classes~\cite{M}, to design differential forms it will be
reasonable to use invariant polynomials on the matrix algebra.
Such polynomial were the trace and the determinant, i.e. forms suggested by
Weyl and Caratheodory respectively. We consider remaining
generators of the ring of invariant symmetrical polynomials from
roots of matrices $\|\frac {\partial S^i}{\partial t^j}\|$.

\section{k-Lagrangian Manifolds}

Let $C$ be a square  $(n \times n)$-matrix.

\begin{definition}.

The minor of a matrix $C$  is called the principal one if the indices
of its rows coincide with indices of its columns.
\end{definition}

\smallskip

\begin{lemma}.

The sum of all principal minors having the order $k$ for a matrix
$C$ is equal to the coefficient of the characteristic polynomial
$P(\lambda)$ standing before the term $\lambda^{n-k}$.
\end{lemma}

\Proof.

It is sufficient only to see on the matrix

$$
C = \left \| c^i_j - \lambda  \delta^i_j\right \|
$$
and to note that the coefficient before $\lambda^{n-k}$ is
obtained by elimination rows and columns of any $k$ diagonal
elements (that leads to principal minors), and subsequent summation
over all such elements.

 $\Box$

Suppose that $n$ action-functions $S^i(t,x), \; (i = 1,...n)$ are
given. To construct an analog of a Hilbert integral, consider the
following closed differential form

\beq
 {\goth S}_k = \sum_{|K|=k} (-1)^{r(K)}
dS^{i_1} ... dS^{i_k}dt^{I \setminus K},
 \label{15}
 \eeq

\noindent where $K = (i_1,...i_k)$.  By $r(K)$ we denote the
number of permutations needed to put members of the exterior
product of differential forms $dt^{i_1} ... dt^{i_k}dt^{I
\setminus K}$ in ascending order. It is easy to calculate that

$$
r(K) = \sum_{i_s \in K}i_s - \frac {k(k+1)}{2}.
$$

The differential form ${\goth S}_k$ corresponds to the coefficient
of $\lambda^{n-k}$ of the characteristic polynomial $P_C(\lambda)$
of the matrix

$$
C = \left \| \frac {d S^i}{d t^j}.
    \right \| = \left \| \frac {\partial S^i}{\partial t^j} +
\frac {\partial S^i}{\partial x^\alpha} \frac {dx^\alpha}{dt^j}
\right \|
$$

Indeed, each of the summands of the matrix ${\goth S}_k$ has the
following Jacobian that corresponds to the principal minor of the
matrix $C$ as a coefficient

$$
 \left | \frac {D[S^K]}{D[t^K]} \right |.
$$

For $k=1$ we have Weyl theory, for $k=n$ Caratheodory theory.

\smallskip

To represent the differential form ${\goth S}_k$ in terms of the
integrand $f$ we need a small regression.

\bigskip

Consider $(n \times n)$-matrix $C = \Phi + \Psi$, where
$\Phi=\diag (\phi^k)$ is a diagonal matrix, $\Psi=\|a^ib_j\|$ is a
rank 1 matrix being the tensor product of a contravariant vector
$a^i$ by a covariant vector $b_j$. We will consider only the case
(needed in what follows) when $\phi^k, a^i$  are differential
forms of the first order, and $b_j$ are elements of the basic
field. Due to the noncommutativity of the exterior product of
differential forms, we will conceive determinant consisted in that
elements in a special, not usual, sense. While expanding the
determinant, we will order factors of each terms in accordance
with numbers of its rows. Taking such rule, a determinant with two
equal column (but not rows!) will be equal to zero. Hence, it is
possible without changing of the determinant to add its columns
(but one cannot add rows!).

\bigskip

\begin{lemma}
The determinant of the matrix $C$ is equal

$$
 \det C = \prod_{i=1}^n\phi^i  +
 \sum_{j=1}^n b_j\phi^1...\phi^{j-1} a^j \phi^{j+1}...\phi^n.
$$

\end{lemma}

\Proof.

Since among the components of the vector $b_j$ there exists at
least one distinct from zero (otherwise $\det C = 0$), we can
set without loss of generality that $b_1 \ne 0$. Preserving the
first column of the matrix $C$ we subtract from each column with
the number $i$ the first one being multiplied by $b_i/b_1$. The
expansion of the matrix obtained relative to the first row gives
the needed formulae.

$\Box$

\smallskip

To find the analog of the Poincare-Cartan form, we rewrite
integrand of~(\ref{1}) as follows

$$
 \frac {1}{{n\choose k}f^{k-1}}\sum_{|K|=k}(-1)^{r(K)}
\sum_{i \in K}\det (diag fdt^i)dt^{I \setminus K} = \frac
{1}{{n\choose k} f^{k-1}}\sum_{|K|=k} \left ( (-1)^{r(K)} \prod_{i
\in K}fdt_i \right ) dt^{I \setminus K}.
$$

Then we add to each factor, standing under the sign of product,
the canonical form for the distribution (\ref{2}):
$\omega^{\alpha} := dx^{\alpha} - \sum_{j=1}^n g^{\alpha}_jdt^j$
being multiplied by the corresponding momentum $p^i_\alpha :=
f_{q^{\alpha}_i}$. The lemma 2 gives the reason to consider as a
 Poincare-Cartan form the following expression

\beq
\Delta = \sum_{|K|=k} \Delta_K,
 \label{16}
\eeq

\noindent where

\beq
\Delta_K
 = \frac {(-1)^{r(K)}}{{n\choose k}f^{k-1}} \prod_{i \in K} \left [ fdt^i +
 \frac {n}{k}  \sum_{\alpha = 1}^ \nu p^i_{\alpha}\left ( dx^{\alpha} -
g^{\alpha}_i (t,x)dt^i \right ) \right ]dt^{I \setminus K}
 \label{17}
 \eeq
Here $\prod $ stands for exterior product of differential forms,
and into the function $f$ we substitute the
distribution~(\ref{2}). In parentheses of the formulae~(\ref{17})
only the summand $g^{\alpha}_idt^i$ takes part, because other
summands of the type $g^{\alpha}_j dt^j$ give zero for $j \ne i$,
since the term $dt^j$ meet twice in the exterior product.

The choice of differential form~(\ref{16}) can be substantiated as
follows. Each summand of the differential form $\sum_K
dS^Kdt^{I\setminus K}$ is a simple multivector of an order $k$. We
can write it as an exterior product of one-dimensional forms.
Using the invariance, we select in each factor, entering in
$dS^K$, independent summands and take it as $dt^i$. The rest can
be justify as follows. The differential form $\Delta$ is reduced
to $fdt$ on the integral manifold of the distribution (\ref{2}).
Besides, the derivatives of integrands of both functionals
relative to $q^{\alpha}_i$ coincide. These facts are provide by
the choice of numerical coefficients ${n\choose k}$ and $n/k$. It
is necessary to speak about minimum of the Weierstrass function
that will be built using the differential form~(\ref{16}) (see below
the section "Weierstrass function").

\begin{note}
The exchange of the determinant by this product is equivalent, in
essence, to the restriction of a symmetrical multilinear
$Ad$-invariant form, defined on a Lie algebra, to the Cartan
subalgebra. This restriction uniquely defines this form~\cite{N}.
\end{note}

Suppose that the distribution (\ref{2}) is integrable. Then the
manifold ${\goth M} \subset J_1(\xi)$, $\dim {\goth M} = \nu + n$
defined by equations~(\ref{2}) is fibred by $n$-dimensional
fibers. Since ${\goth M}$ (due to its definition) has one-to-one
projection onto the space $\xi$, the foliation on ${\goth M}$ induces a
foliation of $\xi$ by $n$-dimensional fibers $\Phi$.

\smallskip

\begin{definition}.

A manifold ${\goth M} \subset J_1(\xi)$, $\dim {\goth M} = \nu +
n$ is called $k$-Lagrangian if the restriction of the differential
form $\Delta$ to ${\goth M}$ is the closed form.

\end{definition}

\smallskip

The differential form $\Delta$ lead to the natural generalization
of the Legenre transform. It corresponds to describing of the
function $f$ from  matrix variable as an envelope of the family of
all tangent planes to the surface $f: \R^{n\nu} \to \R$ considered
as a function of exterior $k$-power of its argument $\frac
{Dx}{Dt}$.

\smallskip

\begin{theorem}

Let the distribution (\ref{2}) be integrable, and the manifold
${\goth M}$ be $k$-Lagrangian. Then the fibers $\Phi \subset \xi$
are solutions to the Euler equations

\beq
 \frac {d}{dt^i} \left ( f_{q^{\alpha}_i} \right ) +
f_{x^\alpha} = 0.
 \label{18}
 \eeq

\end{theorem}

\Proof.

Calculate the differential of the form
$$
\Delta_K = \frac {1}{{n\choose k}}fdt + \frac {1}{{n-1\choose
k-1}} \sum_{{\alpha}=1}^\nu\sum_{s=1}^k [f_{q^{\alpha}_{i_s}}
dt^{i_1}...dt^{i_{s-1}}(dx^{\alpha}-g^{\alpha}_{i_s}dt^{i_s})
dt^{i_{s+1}}...dt^{i_k}]dt^{I \setminus K}.
$$

Summands which contain more than one form $\omega^\alpha$ are
omitted in this formula because after differentiation and
subsequent substitution the distribution (\ref{2}) they give zero.
While differentiating, we consider only summands of the type
$dx^\lambda dt^I$. The result will be presented as four groups.

$$
\begin{array}{c}
\phi_1=\frac {1}{{n\choose k}}\frac {\partial f}
{\partial x^\lambda}dx^\lambda dt^I,\\

\phi_2=\frac {1}{{n\choose k}}f_{q^{\alpha}_j}
\frac {\partial g^{\alpha}_j}{\partial x^\lambda}dx^\lambda dt^I, \\

\phi_3=\sum_{s=1}^k \frac {1}{{n-1\choose k-1}}
f_{q^{\alpha}_{i_s}} \left ( - \frac {\partial
g^{\alpha}_{i_s}}{\partial x^\lambda} \right )
dx^\lambda dt^I \\

\phi_4=\sum_{s=1}^k \frac {1}{{n-1\choose k-1}}\frac {d}{dt^{i_s}}
f_{q^\lambda_{i_s}}dt^{i_s}dt^{i_1}...dt^{i_{s-1}}dx^\lambda
dt^{i_{s+1}}...
dt^{i_k}dt^{I \setminus K}. \\
\end{array}
$$

The sums $\phi_1$ and $\phi_2$ are generated from differentiation
with respect to $x^\lambda$ of the summand $\frac {1}{{n\choose
k}}fdt^I$. The sum $\phi_3$ is generated from differentiation
with respect to $x^\lambda$ of the factor $\omega^{\alpha_s}$. In the
sum $\phi_4$ stand the full derivatives
with respect to $t^{i_s}$. From $\omega^{\alpha}$ is retained only
$dx^\lambda$. To order
differentials in $\phi_4$, it is necessary to permute $dt^{i_s}$
and $dx^\lambda$, so the sum $\phi_4$ changes its sign.

Further we have to sum over $K$ the expression obtained. The
expression $\sum_K \phi_1$ as well as $\sum_K \phi_2$ consist of
equal summands, and its number is equal to the number of groups
$K$. Hence, the coefficients ${n\choose k}$ are cancelled. Each
summand in the expression $\sum_I \phi_3$ as well as in $\sum_I
\phi_4$ is met as many times as there are groups that it contain.
The number of such groups equal ${n-1\choose k-1}$. Hence, the
coefficients $\frac {1}{{n-1\choose k-1}}$ are cancelled. After
this reduction  $\sum_K \phi_2 + \sum_K \phi_3 = 0$. The remaining
summands $\sum_K \phi_1 + \sum_K \phi_4$ give

$$
 \left ( -\frac {d}{dt^i}
\left ( f_{q^{\lambda}_i}
 \right ) + f_{x^\lambda}
\right ) dx^\lambda dt^I = 0.
$$

       $\Box$

\smallskip

\begin{note}
For $n=1$ the set of solutions to Euler equation is
finite-dimensional. For $n>1$ there are functional freedom to
choose solutions. These solutions could be combined in the
"Lagrangian manifolds" variously. This is the reason why there are
many different kinds of fields of extremals for the case of
multiple integral.
\end{note}

\smallskip

\section{Weierstrass function}

Having invariant integral of the Hilbert type corresponding to the
form~(\ref{16}) one can build the analog of the Weierstrass
function for the distribution of the slope field
$g^{\alpha}_i(t,x)$ \cite{I}. To have nonnegative Weierstrass
function for all values of derivatives it is necessary that: first,
the value of the integrand corresponding to Hilbert integral on
manifolds imbedded into the slope field $g^{\alpha}_i(t,x)$ was
equal to the value of $f$, and second, its derivatives
with respect to
$Dx/Dt$ were equal to that of the function $f$. Canonical forms
$\omega^{\alpha}$ vanish on manifolds imbedded into the slope
field in question. Hence, calculating differential one can ignore
summands containing products of form $\omega^{\alpha}_j = \frac
{dx^\alpha}{dt^j} - g^{\alpha}_j(t,x)$. Let us denote by $\hat f$
the value of the function $f$ after substitution of the slope field
~(\ref{2}).

\beq
 {\cal {E}}\left ( t,x,\frac {dx}{dt},g \right ) = f\left (
 t,x,\frac {dx}{dt}\right ) - \sum_{|K|=k}
\frac {1}{{n\choose k} \hat f^{k-1}} \det \|\hat f\I
+\frac {n}{k}\hat p^i_\alpha (\omega^\alpha_j) \|.
 \label{19}
 \eeq

It is evident that the coefficient $\frac {1}{{n\choose k}\hat
f^{k-1}}$ standing before the sum provides the condition $\hat f
=f$. Coefficients standing inside of determinants provide
coincidence of derivatives. Indeed, let us expand the determinant
in~(\ref{19}). Summands containing $\omega^{\alpha}_j = \frac
{dx^\alpha}{dt^j} - g^{\alpha}_j(t,x)$ in the first degree arise
if we take only one factor standing on diagonal. The derivative
$\frac {dx^\mu}{dt^m}$ stands on the $m$ place of diagonal, and
its coefficient is $\frac {n}{k}q^m_\mu$. It does not depend on
$K$. This coefficient meets as many times as there are multiindices
$K$ containing index $m$, that is ${n-1\choose k-1}$  times. The
total coefficient will be equal to ${n-1\choose k-1}/ {n\choose k}
= \frac {k}{n}$. This gives the coincidence of derivatives of the
function $f$ and of the subtrahend in the formula~(\ref{19}).
Hence, the following relations will be valid

\beq
 {\cal {E}}(t,x,g,g) = 0; \quad \frac {\partial {\cal
{E}}}{\partial (\frac {dx^\alpha}{dt^i})}(t,x,g,g) = 0.
 \label{20}
\eeq

\smallskip

\begin{definition}.
A slope field $g^{\alpha}_i (t,x)$ is called geodesic for the
differential form ${\goth S}_k$ in a domain ${\cal V} \subset \xi$
if the minimum of the Weierstrass function~(\ref{19}) is reached
at each point $(t,x) \in {\cal V}$ for $\frac {dx}{dt} = g$.
\end{definition}

\begin{note}
Let us recall that the Legendre transform in case of simple integrals
and the Weyl transform in case of multiple integrals describe a function
$f(\cdot,\cdot,Dx/Dt)$ with the help of support planes to the graph
$f: \R^{n\nu} \to \R$. By contrast, the transforms ${\cal Z}_k$,
defined below in section 6, will  describe $f$ with the help of
multilinear (relative to variables $\partial x^\alpha /\partial t^i$) support
manifolds to the graph $f: \R^{n\nu} \to \R$. These support
manifolds can be regarded as planes in the space of multivectors
of dimension $k$.
\end{note}

\smallskip

\begin{theorem}.

If a manifold $\hat x(\cdot)$ is imbedded into the geodesic field
$g^{\alpha}_i (t,x)$ in the domain ${\cal V}$ then the functional
in question reaches on this manifold the minimal value relative to
any manifold with the same boundary lying in the domain ${\cal V}$.

In addition, the local minimum relative to $dx/dt$ gives the
sufficient condition for the weak local minimum while the global
minimum gives the sufficient condition for the strong minimum in
${\cal V}$.
\end{theorem}

\Proof.

The proof follows by the standard way from the invariance of the
Hilbert integral and from the positivity of the Weierstrass
function in the given domain.

       $\Box$

\smallskip

Note that in this theorem we do not suppose that the slope field
$g^{\alpha}_i (t,x)$ is integrable. Nevertheless, the imbedding
into a geodesic field gives minimum for one individual manifold.
If the field is integrable then there exists a manifold with the
slope defined by the field which passes through each point of the
domain ${\cal V}$. Due to the theorem 2 it gives the minimum, and,
a fortiori, it is the extremal. Hence, in the integrable case, the
manifold $\hat x(\cdot)$ is imbedded into a field of extremals.
Through each point of the domain ${\cal V}$ passes one and only one
$n$-dimensional extremal giving minimum to the functional. The
$n$-parametric family of $\nu$-dimensional level
surfaces of action-functions $S^i = Const$ transversally
cross these extremals.  (Transversality
conditions will be found below in section 5.) The German classical literature
on variational calculus awards to the described geometrical object
the name "Perfect or complete picture" (eine vollst\"andige
Figure).

\bigskip

\section{Condition of solvability}

Let us find conditions for minimum of function ${\cal {E}}$ at a
point $g$. The first derivative at a point $g$ must be
zero. We differentiate $\cal E$ by $dx^\lambda /dt^m$. This argument
enters only into the $m$-column of each determinant. The derivative
will be equal to the same determinant in which the $m$-column is
changed by

$$
\frac {n}{k}\hat p_{\lambda}^i.
$$

By expanding the determinant relative to this column one obtains

\beq
 \hat p^l_{\lambda} - \frac {1}{{n\choose k}\hat f^{k-1}}
\sum_{|K|=k, K \ni l} \sum_{i \in K}\left (  \adj_i^l(C(K))
 \frac {n}{k}\hat p_{\lambda}^i \right ) = 0.
 \label{21}
 \eeq

Denote by $C(K)$ the matrix corresponding to the principal minor
with the index $K$ in formula~(\ref {19}). Denote by $\adj_i^l(C(K))$
the adjunct (the cofactor) of the element $(i,l)$ of the matrix
 $C(K)$. Denote by $\adj_{ij}^{lm}(C(K))$ the adjunct of the minor
standing in the columns ($l,m$) and in the rows ($i,j$) of the matrix
 $C(K)$.

The second derivative of the function ${\cal {E}}$ at a point $g$
is obtained by the exchange two columns by coefficients of the
corresponding $dx/dt$. The expansion of these determinant relative
to pair of columns ($l,m$) gives the following condition of minimum:
The quadratic form with the coefficients

\beq
 \frac {\partial^2 f}{\partial q^{\lambda}_l \partial q^{\mu}_m}
 - \frac {1}{{n\choose k}\bar f^{k-1}} \sum_{|K|=k, K\ni
(lm)}\sum_{(i,j) \in K}\frac {n^2}{k^2}\left ( \adj_{ij}^{lm}
(C(K))\right )  [\hat p_{\lambda}^i \hat p_{\mu}^j - \hat
p_{\lambda}^j \hat p_{\mu}^i]
 \label{22}
 \eeq

\noindent on the space $(n\nu \times n\nu)$-matrix must be
nonnegative. Let us rearrange coefficients of quadratic form
~(\ref{22}) which adds to the first summand

 \beq
\frac {\partial^2 f}{\partial q^{\lambda}_l \partial q^{\mu}_m}.
 \label{23}
 \eeq

If the minor $\adj^{lm}_{ij} C(K)$ is not the principal one
($(l,m) \ne (i,j)$) then, after the substitution $q=g$, it will have
zero row or zero column. If $(l,m)=(i,j)$ then it turns into the
diagonal matrix with the diagonal elements $\hat f$. So, we have
 $\adj^{ij}_{ij} C(K) = \hat f^{k-2}$. Hence, all summands of
coefficients in question appear the same. The number of
these summands is equal to the number of minors $K$ containing
the pair of indices $l,m$. So the total coefficient is equal to
${n-2\choose k-2}$. Hence, we have

\beq
 \frac {\partial^2 {\cal {E}}}{\partial q^{\lambda}_l \partial
q^{\mu}_m}= \frac {\partial^2 \hat f}{\partial q^{\lambda}_l
\partial q^{\mu}_m}- \frac {n(k-1)}{k(n-1)\hat f} (p^l_\lambda
p^m_\mu - p^l_\mu p^m_\lambda) > 0.
 \label{24}
\eeq

The expression $(p^l_\lambda p^m_\mu - p^l_\mu p^m_\lambda)$
defines the skew-symmetric form vanishing on matrices of rank 1.
Hence, its addition do not violates Hadamard-Legendre condition.
The formula~(\ref{24}) gives precisely such skew-symmetrical
forms which has to be added to the integrand to obtain the
invariant integral of one or other degree $k$.

\bigskip

\section{Transversality condition}

\bigskip

Consider moving boundary problems of minimization of the
functional~(\ref{1}). For example we can put that the boundary
$\partial V$ of unknown solution $\hat x(\cdot)$ belongs to a
fixed manifold $\cal X$. It is evident that the solution to this
problem meets the necessary conditions for optimality for the
problem with the fixed boundary. Find additional necessary
conditions caused by the possibility to vary the boundary.

Suppose that each point of the manifold $\hat x(\cdot)$ moves in
space $(t,x)$ along trajectories of a vector field
$T^i(t,x), X^\alpha (t,x)$ that tangent to the manifold $\cal X$.
Denote by $\theta$ the time of translation. The solution of the system

\beq
\left \{
\begin{array}{r}
\dot t = T(t,x)   \\

\dot x = X(t,x) \\
\end{array}
\right.
 \label{25}
\eeq
 \noindent  with initial conditions $(t_0,x_0) \in \hat
x(\cdot)$ will be denoted by ${\cal T}(\theta ;t_0,x_0), {\goth
X}(\theta ;t_0,x_0)$. By $h(t)$ we denote the derivative of $x$
with respect to parameter $\theta$ for given $t$.

$$
h(t) = \frac {\partial \goth X}{\partial \theta}(0;t,\hat x(t)).
$$

After substitution ${\cal T}(\theta),{\goth X}(\theta)$ we obtain
the function ${\cal F}(\theta)$. Find

\beq
  \frac {d}{d\theta}{\cal F}(0) = \int_V \left ( \hat
f_{x^\alpha}h^\alpha + \hat f_{\frac {\partial
x^{\alpha}}{\partial t^i}}\frac {\partial h^{\alpha}}{\partial
t^i} \right ) dt^I + \int_{\partial V}\hat f T_u dS.
\label{26}
 \eeq
\noindent Here $T_u$ is the projection of the vector $T$ on $u$
(on the normal to the boundary $\partial V$), and $dS$ is the
element of the volume of the boundary. Integrate by parts the
summand including derivatives of $h$. Under the integral on $V$ we get
the left hand side of the Euler equation on the manifold $\hat x(\cdot)$
that gives zero. Under the integral on $\partial V$ it is added the summand
$\hat f_{\frac {\partial x^{\alpha}}{\partial t^i}}
h^\alpha (-1)^{i-1}dt^{I \setminus i}$. Let us express $h^\alpha$
at points $\partial V$ through components of the vector field
$(T,X)$.

Differentiation of the identity ${\goth X}(\theta) = x({\cal T}(\theta),
\theta)$ gives

$$
\frac {\partial x^\alpha}{\partial \theta} + \frac {\partial \hat
x^\alpha}{\partial t^j}\frac {\partial t^j}{\partial \theta} =
X^\alpha
$$
\noindent or

 \beq
h^\alpha = X^\alpha - \frac {\partial \hat x^\alpha}{\partial
t^j}T^j.
 \label{27}
 \eeq

Substitute~(\ref{27}) into ~(\ref{26}) and take into account that
$T_\nu dS = T^i(-1)^{i-1}dt^{I \setminus i}$. We obtain

 \beq
\frac {d}{d\theta}{\cal F}(0) = \int_{\partial V} \left [ \hat
f_{\frac {\partial x^{\alpha}}{\partial t^i}}X^\alpha + \left (
\hat f \delta^i_j - \hat f_{\frac {\partial x^{\alpha}}{\partial
t^i}}\frac {\partial \hat x^{\alpha}}{\partial t^j} \right ) T^j
\right ] (-1)^{i-1}dt^{I \setminus i}
 \label{28}
 \eeq

\begin{theorem}

The necessary condition for optimality for moving boundary problem
is

$$
\hat f_{\frac {\partial x^{\alpha}}{\partial t^i}}X^\alpha + \left
( \hat f\delta^i_j - \hat f_{\frac {\partial x^{\alpha}}{\partial
t^i}}\frac {\partial \hat x^{\alpha}}{\partial t^j} \right )T^j =
0
$$
\noindent for any vector field $T(t,x),X(t,x)$ which tangent to the
manifold $\cal X$.
\end{theorem}

\Proof.

Let the theorem be violated for a vector $T(t_0,x_0),X(t_0,x_0)$
which tangent to $\cal X$ at a point $(t_0,x_0)$. Extend
this vector to a smooth field which tangent to $\cal X$ and which
is nonzero only in the sufficiently small neighborhood of the point
$(t_0,x_0)$. For the variation corresponding to the shift along this
vector field we get

$$
\frac {d}{d\theta}{\cal F}(0) < 0,
$$
\noindent that contradicts to minimality of $\hat x$.

$\Box$

The transversality conditions give $n$ equations on the vector
$(T,X)$. The number of conditions is equal to the number of
action-functions. Clarify the meaning of these conditions in case
$k=n$ (the case of Caratheodory). We have

$$
\int_Vfdt^I = \int_V \det \left \| \frac {\partial S^i}{\partial
t^j}  + \frac {\partial S^i}{\partial x^\alpha} \frac {\partial
x^\alpha} {\partial t^j}\right \| dt^I.
$$

Denote the matrix standing under the sign of determinant by $F$.
Denote by $a^r_s$ the matrix appearing in the transversality
condition

$$
a^r_s = f\delta^r_s - f\sum_{i=1}^n (F^{-1})^r_i \frac {\partial
S^i}{\partial x^\alpha} \frac {\partial x^\alpha}{\partial t^s}
$$

The transversality condition is

$$
a^r_s T^s + p^r_\alpha X^\alpha = 0.
$$
\noindent In our case

$$
 \left [ f\delta^r_s - f\sum_{i=1}^n (F^{-1})^r_i \frac {\partial
 S^i}{\partial x^\alpha} \frac {\partial x^\alpha}{\partial t^s}
\right ] T^s +  f \sum_{i=1}^n (F^{-1})^r_i \frac {\partial
S^i}{\partial x^\alpha}
 X^\alpha = 0.
$$

Divide by $f$ and multiply from the left by the matrix $F^m_r$:

$$
\left [ F^m_s - \frac {\partial S^m}{\partial x^\alpha} \frac
{\partial x^\alpha }{\partial t^s} \right ] T^s + \frac {\partial
S^m}{\partial x^\alpha} X^\alpha = 0.
$$

Substitute the value of $F$

$$
\frac {\partial S^m}{\partial t^s}T^s + \frac {\partial
S^m}{\partial x^\alpha}X^\alpha = 0.
$$

Consequently, the vector $(T,X)$ is lying on the intersection of
manifolds $S^m = Const.$

\bigskip

\section{Generalization of Legendre-Weyl-Caratheodory transforms}

\smallskip

We fix value $(t,x)$ as for the classical Legendre transform.
These variables will not appears in subsequent formulas, and by $f(q)$
we always understand $f(t,x,q)$. Recall that

$$
q^\alpha _i = \frac {\partial x^\alpha}{\partial t^i}
$$
\noindent are the main arguments of the function $f$, and

$$
p_\alpha ^i = \frac {\partial f}{\partial q^\alpha _i}
$$
are the corresponding moments. Denote by $\Lambda^kR$ exterior $k$-power
of the matrix $R$. It is the matrix consisted from minors of order $k$
of the matrix $R$. The function that corresponds to the invariant integral~
(\ref{15}) will be written in the form

\beq
 \tilde {\Delta} = \frac {1}{{n\choose k}\hat f^{k-1}}\tr \Lambda^k\left \|
 \hat f\delta^i_j + \frac {n}{k} p_\alpha ^i q^\alpha _j.
  \right \|.
 \label{29}
 \eeq

In view of lemma 1 it is natural to relate the expression $\Delta$
with the function

 \beq
\tilde {\goth S} = \tr \Lambda^k
 \left \| \frac {\partial S^i}{\partial
t^j} + \sum_{\alpha=1}^\nu \frac {\partial S^i}{\partial x^\alpha}
\frac {\partial x^\alpha}{\partial t^j} \right \|.
 \label{30}
\eeq

The conjugate variables (with the variable $q^\alpha _i$) which will
play the part of new independent variable in our analog of Legenre
transform must be tied with $\partial S^i/\partial x^\alpha$. Keeping
this in mind we introduce the matrix

\beq
 (q^*)_m^\mu (t,x,q) = \left (\frac {\partial \tilde {\Delta}}
{\partial (q_m^\mu)} \right )^t ,
 \label{31}
\eeq
The differentiation of~(\ref{29}) gives

$$
(q^*)_m^\mu =  \frac {1}{{n\choose k}\hat f^{k-1}}\sum_{K \ni (m,l)}
\left ( \adj ^m_l
 \left \| \hat f\delta^i_j + \frac {n}{k}p_\alpha ^iq^\alpha _j
 \right \|_K  \right )^t \frac {n}{k}[(p)^t]_l^\mu.
$$
The subscript $K$ means that entries of the matrix have indices belonging to $K$.
Using the auxiliary matrix

 \beq
[A^{-1}]_m^l =  \frac {1}{{n\choose k}f^{k-1}}\sum_{K \ni (m,l)} \left (
 \adj^m_l \left \| f
 \delta^i_j + \frac {n}{k}p_\alpha ^iq^\alpha _j
 \right \| \right )_K \frac {n}{k}
 \label{32}
 \eeq
we write $(q^*)^\mu _m$ in the form

\beq
 (q^*)_m^\mu = [p^t]_l^\mu [A^{-1}]_m^l.
 \label{33}
 \eeq

Denote the canonical function corresponding to $f$ by
\beq
H = f^*(q^*) = \gamma (\tilde {\Delta} - f).
\label{34}
 \eeq
The definition of this function will be given later in the formula (38).

Let us express the differential form $\Omega$ related to the function
$\tilde {\Delta}$ in canonical coordinate. We carry out of brackets
the function $f^*(q^*)$. Each simple multi-vector entering
in the sum  can be represented as a product of $k$
differential forms of the first order. Select in each form the summand
$dt^i, \; (i \in K)$. The number of these simple multivectors is
${n\choose k}$ and we cancel $\Omega$ by this factor. The coefficient
of $q_i^\alpha$ has to
be equal to $\hat f_{q_i^\alpha}$. Consequently, coefficients of
$dx^\alpha$ in factors of each monomial include ${n-1\choose k-1}$
--- the number of minors $K$ containing the index $i$. In view of our
cancellation by ${n\choose k}$ this factor is $\frac {n}{k}$. As a
result, under the integral appears multilinear function
(relative to variables $q_i^\alpha$). Its derivative with respect to
$q_m^\mu$ gives $(q^*)^m_\mu$. Consequently,
this function coincides with the function $\Delta$. We arrive to the following
differential form

\beq
 \Omega =  H \sum_{|K|=k} (-1)^{r(K)} \left (
 \prod\limits_{i \in K}(dt^i + \sum_{\alpha =1}^\nu \frac {n[(q^*)^t]^i_\alpha}
 {kH}dx^\alpha) dt^{I\setminus K} \right ).
 \label{35}
 \eeq

Normalize variables by putting

\beq
 Q^i_\alpha = \frac {n[(q^*)^t]^i_\alpha}{kH}.
 \label{36}
 \eeq
The formula~(\ref{35}) takes the form

\beq
 \Omega =  H \sum_{|K|=k} (-1)^{r(K)}\left (
 \prod\limits_{i \in K}(dt^i + \sum_{\alpha =1}^\nu Q^i_\alpha
 dx^\alpha) dt^{I\setminus K} \right ).
 \label{37}
 \eeq

The condition of minimum of the Weierstrass function~(\ref{21}) can
be written as the following equation
$p^i_\alpha = \partial \tilde {\Delta} /\partial q^\alpha _i$, or $q^* = p^tA^{-1}$.
In view of~(\ref{24}) the equation~(\ref{31}) can be resolved relative
to $q$. It gives the dependence $q = \phi (q^*)$.

Define the transform ${\cal Z}_k$ following the Caratheodory"s approach.
Consider quadruples  $\{f, \phi , q_i^\alpha , p^i_\alpha\}$, which satisfy the
condition
$$
f + \phi = \tilde {\Delta}.
$$

Introduce the auxiliary matrix $A$ by using~(\ref{32}). In view of~(\ref{24})
the matrix $A$ on extremals has to be nonsingular. The formula~(\ref{33}) gives

$$
 \frac {\partial \tilde {\Delta}}{\partial q_m^\mu} = [p^t]_i^\mu
(A^{-1})_m^i.
$$

Denote by $\gamma$ the scalar coefficient

$$
\gamma = \frac {\hat f^{k-2}}{\det A}.
$$

Variables that are images of quadruples under the transform ${\cal Z}_k$
will be marked by superscript *. Define ${\cal Z}_k$ by formulas

\beq
 f^* = \gamma \phi ; \; \phi^* = \gamma f; \; (p^*) = \gamma
Aq^t; \; (q^*)=p^tA^{-1}.
 \label{38}
  \eeq
The transformed function will be $f^*(q^*)$.

\begin{theorem}
The transform ${\cal Z}_k: f(q) \mapsto f^*(q^*)$ is birational and involutory.
\end{theorem}

\Proof.

The fact that ${\cal Z}_k$ is birational is evident.

The proof of involutority.

$$
q^*p^* = \gamma (p^tA^{-1})(Aq^t) =\gamma p^tq^t.
$$

Since the transposition does not changes determinants we have

\beq
 \tilde \Delta^* = \gamma  \tilde \Delta .
 \label{39}
\eeq

One obtains by multiplication of the formula $f + \phi = \tilde {\Delta}$
by $\gamma$
$$
f^* + \phi^* = \tilde {\Delta}^*.
$$

The auxiliary matrix $A^{-1}$ was defined through differentiation of
the function~(\ref{29}) Therefore the matrix $(A^{-1})^*$ will be defined
through differentiation of the matrix

$$
 \frac {1}{{n\choose k}(f^*)^{k-1}}\tr \Lambda^k\left \|
 (f^*)\delta^i_j + \frac {n}{k} (p^*)_\alpha ^i (q^*)^\alpha _j.
  \right \|.
$$

In view of the formulas~(\ref{38}), ~(\ref{39}), this matrix equals

$$
 \frac {\gamma}{{n\choose k}f^{k-1}}\tr \Lambda^k\left \|
 f\delta^i_j + \frac {n}{k} p_\alpha ^i q^\alpha _j.
  \right \|^t.
$$

Hence $A^* = \gamma A^t$, and
$$
\gamma^* = \frac {(f^*)^{k-2}}{\det A^*} =
 \frac {\gamma^{k-2} (f)^{k-2}}{\gamma^k\det A}  = \frac {1}{\gamma}.
$$

The formulas obtained allow us to find the second iteration of
the transpose ${\cal Z}_k$:

\beq
\begin{array}{c}
 f^{**} =\gamma^*(\phi^*) = f; \; \phi^{**} = \gamma^*
f^* = \phi; \\

 p^{**}=\gamma^*A^*(q^*)^t = \gamma^*\gamma A^t(A^{-1})^tp = p; \\
  q^{**}=(p^*)^t(A^*)^{-1} = \gamma (qA^t)^t\frac {1}{\gamma}(A^{-1})^t = q.
\end{array}
 \label{40}
\eeq

The theorem is proved.

 $\Box$

\smallskip

In full agreement with the fact that ${\cal Z}_k$ is involutory
and tangential transform we have that the form $\Omega$ is
obtained from the differential form corresponding to the
function~(\ref{29}) by exchange all the variables by its dual
(in the sense of the transform ${\cal Z}_k$).

\begin{note}
There is an intimacy between the transforms ${\cal Z}_k$ and
the condition of polyconvexity of function $f$. Indeed, if the
function $f$ admits a convex extension to the space $V^k(Dx/Dt)$
then it can be described as an envelope of support planes in
the space of multivectors (see the note 3 in the section 3).
Hence, after the addition of corresponding expression of the type~
(\ref{24}), we can apply to $f$ the transform ${\cal Z}_k$.
\end{note}

\smallskip

\section{Formulas for action-functions}

\smallskip

Let us compare two expression of the integrand of invariant integral.
The first is in terms of action-function

 \beq
\goth S = \sum_{|K|=k}\det \left \| \frac {\partial S^i}{\partial
t^j} + \sum_{\alpha=1}^\nu \frac {\partial S^i}{\partial x^\alpha}
\frac {\partial x^\alpha}{\partial t^j} \right \|.
 \label{41}
\eeq

The second is in terms of initial integrand $f$ being written as
integrand of differential form $\Omega$~(\ref{38}).

\beq
  H\sum_{|K|=k} \left (
 \prod\limits_{i \in K}(1 + \sum_{\alpha =1}^\nu Q^i_\alpha
 \frac {\partial x^\alpha}{\partial t^i}) \right ).
 \label{42}
 \eeq

We expand each of determinants of the formula~(\ref{41}) in
the sum of its columns and gather terms with the determinants
of the same order. To write the corresponding formulas we put
that $J \subset K$ is a multiindex, and
$\Xi \subset \{1,2,...\nu\}$ is a multiindex of the same order
as the order of $J$.

\beq
 \goth S = \sum_{|K|=k} \left ( \left | \frac {{\cal D}[S^K]}{{\cal
 D}[t^K]} \right | + \sum_{s=1}^{\min (k,\nu )}\sum_{|J|=|{\Xi}|=s}\left (
\left |  \frac {{\cal D}[S^{K\setminus J}]}{{\cal D}[t^{K\setminus
J}]}\right | \left | \frac {{\cal D}[S^J]}{{\cal D}[x^\Xi ]}
\right | \left | \frac {{\cal D}[x^\Xi]}{{\cal
 D}[t^J]} \right | \right ) \right ) .
  \label{43}
 \eeq

Multiply brackets under the sign of the product in the formula~(\ref{36}),
taking into account that the Jacobian $\frac {{\cal D}[x^\Xi]}{{\cal  D}[t^J]}$
is equal to the coefficient arising by expressing of $dx^\Xi$ through
$dt$. Comparing~(\ref{38}) with~(\ref{43}) we obtain the series of formulas

$$
\sum_{|K|=k} \left | \frac {{\cal D}[S^K]}{{\cal D}[t^K]} \right |
= H.
$$

$$
\sum_{K\ni i} \left | \frac {{\cal D}[S^K]}{{\cal D}[t^{K\setminus
i}x^\alpha]} \right |  = H \sum_{K\ni i} Q^i_\alpha =
 {n-1\choose k-1} (q^*)^i_\alpha.
$$

Finally, to write the general formula we make denotation more
concrete. We fix the sets $J = \{j_1,...j_s\}$ and $\Xi =
\{\alpha_1,...\alpha_s\}$. Then

$$
\sum_{K\supset J} \left | \frac {{\cal D}[S^K]}{{\cal
D}[t^{K\setminus J}x^\Xi ]} \right |  = H \sum_{K\supset J} \det
\| Q^J_\Xi\| = {n-s\choose k-s}H^{1-s}
 \det \|(q^*)^J_\Xi \|.
$$

Expressions $\frac {{\cal D}[S^K]}{{\cal D}[t^{K\setminus J}x^\Xi ]}$
are the Pl\"ucker coordinates of $k$-multivectors composed from
gradients of the action-functions $S^i$ (relative to both dependent
and independent arguments $(x^\alpha ,t^i)$). Summarize: The
canonical variables $H$ and $Q$ are sums of the Pl\"ucker
coordinates of gradients of the action-functions in the bundle
$\xi$. It is the usual gradient-vector in the Weyl construction,
and it is only one multivector of maximum dimension $n$ --- the
determinant --- in the construction of Caratheodory.

Now we find the corresponding expression using initial function $f$.
Recall that

\beq
  \Delta = \sum_{|K|=k} \frac {1}{{n\choose k}f^{k-1}} \det
\left \| f\delta^i_j + \frac {n}{k}\sum_{\alpha=1}^\nu (p^i_\alpha
(q^\alpha _j -g^\alpha _j)) \right \|.
 \label{44}
\eeq

With the determinant in the formula~(\ref{44}) we shall carry out
the same operation as with the function~(\ref{43}). We expand each
determinant into the sum of its columns and gather determinants of
the same order together.

We calculate coefficients of the corresponding minors.
Coefficient of $f$ under the sign of the sum over $K$ is equal ${n\choose k}$,
and this sum include the same number of identical summands. So, the
total coefficient equals to $1$. Coefficient of $f^0p^i_\alpha
\frac {\partial x^\alpha}{\partial t^i}$ under the sign of sum over
$K$ is equal to ${n-1\choose k-1}$ and the sum over $K$ includes the same
number of identical summands, since it is just the number of $k$-sets
including the index $i$. So, the total coefficient equals to $1$ too.
We have  ${n-2\choose k-2}$ summands including the pair of indices
$(i,j)$, and the coefficient of $f^{-1}\frac {{\cal D}[x^{\alpha
,\beta}]}{{\cal D}[t^{i,j}]}$ is equal to $\frac {k-1}{n-1}$. The
coefficient of the general member $f^{-s+1}\frac {{\cal D}[x^\Xi ]}{{\cal
D}[t^{J}]}$ is equal to ${n-s\choose k-s}/{n-1\choose k-1}$.

The minor composed from $\partial f /\partial p^i_\alpha$, where
$\alpha \in \Xi ,\; i \in J$ will be denoted by

$$
\frac {{\cal D}[f]}{{\cal D}[p_\Xi^J]}.
$$

We have the formula

\beq
 \Delta =  f + p^i_\alpha \frac {\partial x^\alpha}{\partial t^i}
+ \sum_{s=2}^{\min (k,\nu)}\sum_{|J|=|\Xi |=s}\left ( f^{-s+1}
\frac {{n-s\choose k-s}}{{n-1\choose k-1}}\left | \frac {{\cal
D}[f]}{{\cal D}[p_\Xi ^J]}\right | \left | \frac
 {{\cal D}[x^\Xi]}{{\cal D}[t^J]}\right |  \right )
  \label{45}
 \eeq

If we equate coefficients of $\frac {{\cal D}[x^\Xi ]}{{\cal
D}[t^J]}$ in formulas~(\ref{44}) and~(\ref{45}) we obtain
equations for action-functions in terms of $f$.

\smallskip

\section{The canonical equations}

\smallskip

It is natural to write the analog of the Jacobi equation in
the canonical coordinates, i.e. in terms of Pl\"ucker cordinates
of multivectors composed from gradients of the action-functions.
To obtain these equations it is sufficient to write the conditions
of closeness of the differential form $\Omega$:

$$
 \Omega =  H \sum_{|K|=k} (-1)^{r(K)}
 \left (
 \prod\limits_{i \in K}(dt^i + \sum_{\alpha =1}^\nu Q^i_\alpha
 dx^\alpha) dt^{I\setminus K} \right ).
$$

We have

\beq
\begin{array}{c}
  d\Omega = \sum_K\left ( \frac {\partial H} {\partial t^m}dt^m +
  \frac {\partial H}{\partial x^\mu}dx^\mu \right )
 \left ( \sum_K \prod (dt^i + Q^i_\alpha dx^\alpha ) \right ) \\

   + H\sum_K  \left  ( \frac {\partial Q^l_\mu}{\partial x^\lambda}dx^\lambda
  + \frac {\partial Q^l_\mu}{\partial t^m}dt^m  \right )
 \prod(dt^i + Q^i_\alpha dx^\alpha) = 0.
 \end{array}
 \label{46}
  \eeq

All the summands of this product are ordered relative
to the order of $dt^i$
but in factors that do not contain the index $i$, i.e.
$i \notin K$, we omit all summands adding to $dt^i$. Besides,
in factors $dt^l + Q^l_\mu dx^\mu$ we omit that functions $Q^l_\mu$
the derivatives of which was already calculated. We equate
coefficients of the form $d\Omega$ for independent products of
differential. The coefficient of $dx^\mu dt^I$ equals

\beq
 \sum_{K\ni m}\left ( \frac {\partial H}{\partial x^\mu} - Q^m_\mu
 \frac {\partial H}{\partial t^m} -
 H\frac {\partial Q^m_\mu}{\partial t^m} \right ) = 0.
\label{47}
 \eeq
The sign "minus" for two last members in the formula~(\ref{47})
is stipulated by necessity to exchange the order of differentials
$dt^m$ and $dx^\mu$.

All summands standing under the sign of the sum are identical,
hence each of it is zero. To write equations in a more compact
form, it is convenient to introduce an operator

$$
L_\mu = \frac {\partial }{\partial x^\mu} - Q^m_\mu
 \frac {\partial }{\partial t^m}.
$$

In terms of the operator $L$ the equation~(\ref{47}) takes the form

\beq
  H \frac {\partial Q^m_\mu}{\partial t^m} =  L_\mu H.
\label{48}
 \eeq

\begin{note}
The abbreviated writing $L_\alpha$ is not an occasional one.
With the heuristic point of view, it is convenient to imagine
$Q^i_\alpha$ as if it were the differential of $t^i$ with respect to
$x^\alpha$ do not calling attention to dimensions. Then the
operator

$$
L_\alpha = \frac {\partial}{\partial x^\alpha} - Q^m_\alpha
 \frac{\partial}{\partial t^m}
$$
would be the full derivative with respect to $x^\alpha$ taking into
account the imaginary dependence $t(x)$.
\end{note}

Rewrite the equation in an another form substituting in it
 $Q^m_\mu = \frac {n[(q^*)^t]^m_\mu}{kH}$. We get

$$
\frac {n}{k}H\frac {\partial [(q^*)^t]^m_\mu}{\partial t^m} - \frac
{n}{k} [(q^*)^t]^m_\mu \frac {\partial H}{\partial t^m} =
 H\frac {\partial H}{\partial x^\mu} -
HQ^m_\mu \frac {\partial H}{\partial t^m}
$$
The second members in both sides of the last equation
mutually annihilate. The equation~(\ref{48}) takes the form

\beq
   \frac {n\partial [(q^*)^t]^m_\mu}{k\partial t^m} =
\frac {\partial H}{\partial x^\mu}
 \label{49}
 \eeq

The definition~(\ref{34}) of the function $H$ and the involutivity
of the trasform ${\cal Z}_k$ give

\beq
  - \frac {\partial H}{\partial [(q^*)^t]^m_\mu} = q_m^\mu =
\frac {\partial x^\mu}{\partial t^m} .
 \label{50}
 \eeq

Combining~(\ref{49}) and~(\ref{50}), we get the system
\beq
 \left \{
\begin{array}{l}
 \frac {n}{k} \frac {\partial [(q^*)^t]^m_\mu}{\partial t^m} =
\frac {\partial H}{\partial x^\mu}  \\

   \frac {\partial x^\mu}{\partial t^m}  =
   - \frac {\partial H}{\partial [(q^*)^t]^m_\mu}.
\end{array}
\right .
 \label{51}
 \eeq

\begin{note}
The Euler equations were obtained by equating to zero the
coefficient of $dx^\lambda dt^I$ of the differential form
$d\Delta$. The formula~(\ref{49}) was also obtained by equating
to zero the same coefficient of the same form but being written
in canonical variables (in terms of Pl\"ucker coordinates of
gradients of the action-functions). Hence, it is natural to
consider the system of equations~(\ref{51}) as a canonical
form of the Euler equation. As in the classical case, the
equation~(\ref{50}) follows from very definition of the functions
$H$ and $q^*$. The shape of the canonical system is the same for
any $k$ but the Hamiltonian as well as the function $q^*$ depends
on the choice of the invariant integral.
\end{note}

\bigskip

\section{Necessary and sufficient conditions \\ for closeness of the form
$\Omega$}.

\smallskip

Calculate the coefficient of $dx^\mu dx^\lambda dt^{I\setminus l}$
at the left hand side of the formula~(\ref{46}).

The summand $\frac {\partial H}{\partial t^m}dt^m$ gives two
members. For one of these members we choose the differential $dx^\mu$
standing on the m-spot, and differential $dx^\lambda$ on the
l-spot. To reorder we need $l$ transpositions. For another one
we choose the differential $dx^\mu$
standing on the l-spot, and differential $dx^\lambda$ on the
m-spot. To reorder we need $(l-1)$ transpositions. These two
members take the form

$$
(-1)^l\frac {\partial H}{\partial t^m}(Q^l_\lambda Q^m_\mu -
Q^l_\mu Q^m_\lambda ).
$$

The summand $\frac {\partial H}{\partial x^\mu}dx^\mu$ gives
two members too. By the same reasoning we verify that its sum
is

$$
(-1)^l\left ( \frac {\partial H}{\partial x^\lambda} Q^l_\mu -
\frac {\partial H}{\partial x^\mu} Q^l_\lambda \right ).
$$

The summand $\frac {\partial Q^l_\mu}{\partial x^\lambda}dx^\lambda$
from the second part of the formula~(\ref{46}) gives the following
pair of members

 $$
(-1)^l H \left ( \frac {\partial Q^l_\mu}{\partial x^\lambda} -
\frac {\partial Q^l_\lambda}{\partial x^\mu} \right ).
 $$

Finally, the summand $\frac {\partial Q^l_\mu}{\partial t^m}dt^m$
gives the fourth pair of members

$$
(-1)^l H \left ( \frac {\partial Q^l_\lambda}{\partial t^m}Q^m_\mu
- \frac {\partial Q^l_\mu}{\partial t^m}Q^m_\lambda \right ).
$$

For the fixed $l$ the expression standing under the sign of
the sum by multiindex $K$ does not depend on $K$. Hence, the
sum of all four written out pairs equals to zero. We get

\beq
 \begin{array}{c}
\frac {\partial H}{\partial t^m}(Q^l_\lambda Q^m_\mu - Q^l_\mu
Q^m_\lambda ) + \left ( \frac {\partial H}{\partial x^\lambda}
Q^l_\mu - \frac {\partial H}{\partial x^\mu} Q^l_\lambda \right )
\\ +  H \left ( \frac {\partial Q^l_\mu}{\partial x^\lambda} -

\frac {\partial Q^l_\lambda}{\partial x^\mu} \right ) + H \left (
\frac {\partial Q^l_\lambda}{\partial t^m}Q^m_\mu - \frac
{\partial Q^l_\mu}{\partial t^m}Q^m_\lambda \right ) = 0.
\end{array}
 \label{52}
\eeq

Using the definition of the operator $L$ and grouping in the due order
the summands we can rewrite the formula~(\ref{52}) in the form

$$
Q^l_\mu L_\lambda (H) + HL_\lambda (Q^l_\mu ) - Q^l_\lambda L_\mu
(H) - HL_\mu (Q^l_\lambda ) = 0.
$$

Since $L$ is the differential operator of the first order we
can use the Leibnitz rule for differentiation of product.
Hence, the last formula turns into

$$
L_\lambda (Q^l_\mu H) - L_\mu (Q^l_\lambda H) = 0.
$$

Taking into account the definition $Q^i_\alpha = \frac {n(q^*)^i_\alpha}{kH}$,
we get the final formula

\beq
 L_\lambda (q^*)^l_\mu = L_\mu (q^*)^l_\lambda ,
 \label{53}
 \eeq
which is the natural generalization of the potential condition
of the vector of momentum.

\bigskip

The subsequent calculations follow the same pattern. We introduce
the results abbreviating a little details of justifications.

Calculate coefficient of $dx^\mu dx^\lambda dx^\rho dt^{I\setminus
\{lm\}}$ at the left hand side of the formula~(\ref{46}).

The first summand $\frac {\partial H}{\partial t^m}dt^m$ standing
before the first sum~(\ref{46}) gives the six members:

\beq
\begin{array}{c}
 (-1)^m \frac {\partial H}{\partial t^m} (Q^m_\rho Q^l_\mu
 Q^r_\lambda - Q^m_\rho Q^l_\lambda Q^r_\mu + Q^m_\lambda Q^l_\rho
 Q^r_\mu \\
 - Q^m_\lambda Q^l_\mu Q^r_\rho + Q^m_\mu Q^l_\lambda
 Q^r_\rho - Q^m_\mu Q^l_\rho Q^r_\lambda ).
\end{array}
 \label{54}
\eeq

The second summand $\frac {\partial H}{\partial x^\mu}dx^\mu$ standing
before the first sum~(\ref{46}) gives also the six members:

\beq
\begin{array}{c}
 (-1)^l \left (
 \frac {\partial H}{\partial x^\mu} (Q^l_\lambda Q^r_\rho
  - Q^l_\rho Q^r_\lambda ) +\frac {\partial H}{\partial x^\lambda}
   (Q^l_\rho Q^r_\mu - Q^l_\mu Q^r_\rho )
  + \frac {\partial H}{\partial x^\rho}  (Q^l_\mu Q^r_\lambda -
   Q^l_\lambda Q^r_\mu )  \right ) .
\end{array}
 \label{55}
\eeq

The summand $H\frac {\partial Q^l_\mu}{\partial
x^\lambda}dx^\lambda$ standing under the sign of the second
sum~(\ref{46}) gives twelve members:

\beq
\begin{array}{c}
 \left (
H\frac {\partial Q^l_\lambda}{\partial x^\mu}Q^r_\rho + H\frac
{\partial Q^r_\rho}{\partial x^\mu}Q^l_\lambda
 \right ) -
 \left (
 H\frac {\partial Q^l_\rho}{\partial x^\mu}Q^r_\lambda +
 H\frac {\partial Q^r_\lambda}{\partial x^\mu}Q^l_\rho
 \right ) + \\
 \left (
  H\frac {\partial Q^l_\rho}{\partial x^\lambda}Q^r_\mu +
  H\frac {\partial Q^r_\mu}{\partial x^\lambda}Q^l_\rho
 \right ) -
 \left (
H\frac {\partial Q^l_\mu}{\partial x^\lambda}Q^r_\rho + H\frac
{\partial Q^r_\rho}{\partial x^\lambda}Q^l_\mu
 \right ) + \\
 \left (
H\frac {\partial Q^l_\mu}{\partial x^\rho}Q^r_\lambda + H\frac
{\partial Q^r_\lambda}{\partial x^\rho}Q^l_\mu
 \right ) -
 \left (
  H\frac {\partial Q^l_\lambda}{\partial x^\rho}Q^r_\mu +
  H\frac {\partial Q^r_\mu}{\partial x^\rho}Q^l_\lambda
 \right ) .
\end{array}
 \label{56}
\eeq

The summand $H\frac {\partial Q^l_\mu}{\partial t^m}dt^m$
standing under the sign of the second sum~(\ref{46}) gives
twelve members:

\beq
\begin{array}{c}
 HQ^m_\lambda \left (
\frac {\partial Q^l_\mu}{\partial t^m}Q^r_\rho + \frac {\partial
Q^r_\rho}{\partial t^m}Q^l_\mu  \right ) -
 HQ^m_\lambda \left (
\frac {\partial Q^l_\rho}{\partial t^m}Q^r_\mu + \frac {\partial
Q^r_\mu}{\partial t^m}Q^l_\rho  \right ) + \\
 HQ^m_\mu \left (
\frac {\partial Q^l_\rho}{\partial t^m}Q^r_\lambda + \frac
{\partial Q^r_\lambda}{\partial t^m}Q^l_\rho  \right ) -
 HQ^m_\mu \left (
  \frac {\partial Q^l_\lambda}{\partial t^m}Q^r_\rho +
\frac {\partial Q^r_\rho}{\partial t^m}Q^l_\lambda
   \right ) + \\
 HQ^m_\rho \left (
\frac {\partial Q^l_\lambda}{\partial t^m}Q^r_\mu + \frac
{\partial Q^r_\mu}{\partial t^m}Q^l_\lambda  \right ) -
 HQ^m_\rho \left (
\frac {\partial Q^l_\mu}{\partial t^m}Q^r_\lambda + \frac
{\partial Q^r_\lambda}{\partial t^m}Q^l_\mu  \right ) .
\end{array}
 \label{57}
\eeq

Gathering together all these summands~(\ref{54})-(\ref{57}), and
using the definition of the operator $L_\alpha$, we get

\beq
 L_\mu [H(Q^l_\lambda Q^r_\rho - Q^l_\rho Q^r_\lambda)] +
L_\lambda [H(Q^l_\rho Q^r_\mu - Q^l_\mu Q^r_\rho)] +
 L_\rho [H(Q^l_\mu Q^r_\lambda - Q^l_\lambda Q^r_\mu)] = 0.
 \label{58}
\eeq

Let us mark that the expression~(\ref{58}) does not depends on
the multi-index $K$. That is the reason why the right hand side
of the formula~(\ref{58}) is zero.

\bigskip

To obtain the general formula we consider two multiindices
$J = \{i_1,...i_s\}\subset \{1,2,...n\}$ of the order $s$ and $\Xi =
\{\alpha_1,...\alpha_{s+1}\}\subset \{1,2,...n\}$ of the order $s+1$.
Let us calculate the coefficient of the member $dx^\Xi
dt^{I\setminus J}$ in the form $d\Omega$. Let us fix the index
$\mu \in \Xi$.

The summand $\frac {\partial H}{\partial x^\mu}dx^\mu$ standing
before the first sum of the formula~(\ref{46}) gives members
consisted of products of $s$ factors, that was obtained from
$Q^i_\alpha dx^\alpha$ where $\alpha \in \Xi \setminus \mu$.
It is necessary to choose the index $i$ from the multiindex $J$,
since otherwise would arise the factor $dt^i$ which does not
include into the product of differentials under consideration.
Each index $\alpha \in \Xi \setminus \mu$ and $i \in J$ must be
taken once and only once, and moreover, a permutation of a pair of indices
$\alpha_i$ and $\alpha_j$ leads to the permutation of the
differentials $dx^\alpha_i$ and $dx^\alpha_j$, and to the changing
of the sign of the coefficient of $\frac {\partial H}{\partial x^\mu}$.
Hence, while ordering the product of the differentials we get the sign
corresponded to the parity of the permutation.
These demands define exactly the determinant

\beq
 \det \|Q^J_{\Xi \setminus \mu} \|,
 \label{59}
\eeq
composed of elements $Q^i_\alpha$, where $i\in J, \; \alpha \in \Xi
\setminus \mu$. To define the sign of this determinant it is
suffice to find the sign of one concrete summands among its
expansion.

\smallskip

Coefficients of members obtained from the summand $\frac {\partial
H}{\partial t^m}dt^m$ standing before the first sum of the
formula~(\ref{46}) give members consisting of products of $s$
factors that is obtained from the summands $Q^m_\alpha dx^\alpha$
with the upper index $m$ (to avoid the arising of the existing
differential $dt^m$). By repeating the previous reasoning we verify
that the sum of all such differentials has as a coefficient

\beq
 \sum_{\mu =1}^\nu \frac {\partial H}{\partial t^m} Q^m_\mu
 \det \|Q^J_{\Xi \setminus \mu} \|,
 \label{60}
\eeq

Joining~(\ref{59}) and (\ref{60}), gives

\beq
\begin{array}{c}
 \sum_{\mu =1}^\nu \left ( \frac {\partial H}{\partial x^\mu} +
  \frac {\partial H}{\partial t^m} Q^m_\mu \right )
 \det \|Q^J_{\Xi \setminus \mu} \| dx^\Xi dt^{I\setminus J} = \\
\sum_{\mu =1}^\nu L_\mu (H)
 \det \|Q^J_{\Xi \setminus \mu} \| dx^\Xi dt^{I\setminus J}.
\end{array}
 \label{61}
\eeq

\smallskip

Let us fix the index $\mu$ and gather the summands that are
obtained from the members $\frac {\partial Q^l_\lambda}{\partial
x^\mu}dx^\mu$ of the second sum in the formula~(\ref{46}).
Coefficients of these summands are adjuncts of the elements
$Q^l_\lambda$ in the matrix
 $\| Q^J_{\Xi \setminus \mu}\|$ with the natural cyclic ordering of
rows and columns. After its multiplication by $\partial
Q^l_\lambda /\partial x^\mu$ and after summation we get the derivative
of the determinant with respect to $x^\mu$. One obtains

\beq
 \sum_{\mu =1}^\nu H \frac {\partial }{\partial x^\mu}
\left ( \det \|Q^J_{\Xi \setminus \mu} \right )  \| dx^\Xi
dt^{I\setminus J}.
 \label{62}
\eeq

\smallskip

Summands obtained from $\frac {\partial Q^l_\mu}{\partial
t^m}dt^m$ standing before the second sum of~(\ref{46}) define
the derivative of the same determinant over $t^m$. This
leads to the following expression

\beq
 \sum_{\mu =1}^\nu H Q^m_\mu \frac {\partial }{\partial t^m}
 \left ( \det \|Q^J_{\Xi \setminus \mu} \| \right )  dx^\Xi
dt^{I\setminus J}.
 \label{63}
\eeq

The addition of~(\ref{62}) and~(\ref{63}) define the action of the
operator $\L_\mu$ on the corresponding determinant. We unify
~(\ref{61}) and (\ref{62}) and get

$$
 \sum_{\mu =1}^\nu  \left (  L_\mu (H)
 \det \|Q^J_{\Xi \setminus \mu} \| dx^\Xi dt^{I\setminus J} +
H L_\mu \left ( \det \|Q^J_{\Xi \setminus \mu} \| \right )
  dx^\Xi dt^{I\setminus J} \right ).
$$

These expressions do not depend on the multiindex $K$, hence,
they equal to zero.

The following theorem was proved

\begin{theorem}

The necessary and sufficient conditions for the differential
form $\Omega$~(\ref{38}) to be closed is

\beq
 \sum_{\mu =1}^\nu  \left (  L_\mu (H)
 \det \|Q^J_{\Xi \setminus \mu} \| dx^\Xi dt^{I\setminus J} +
H L_\mu \left ( \det \|Q^J_{\Xi \setminus \mu} \| \right )
  dx^\Xi dt^{I\setminus J} \right ) =0
 \label{64}
 \eeq
for any choice $\Xi$ and $J$.

\end{theorem}

\bigskip

\section{Connection and Curvature}

\bigskip

To find the connection generated by a field of extremals let
us turn to the canonical system~(\ref{51}).

$$
 \left \{
\begin{array}{l}
 \frac {n}{k} \frac {\partial x^\mu}{\partial t^m}  =
   - \frac {\partial H}{\partial [(q^*)^t]^m_\mu} \\

   \frac {\partial [(q^*)^t]^m_\mu}{\partial t^m} =
\frac {\partial H}{\partial x^\mu}.  \\
\end{array}
\right .
 $$

Consider the  variational equations, i.e. the system for
the derivatives of solutions with respect to a parameter.

\beq
 \frac {n\partial [(q^*)^t]^m_\mu}{k\partial \eta^\sigma} = U^m_{\mu
 \sigma}; \quad  \frac {\partial x^\lambda}{\partial \eta^\sigma}
 = V^\lambda _{\sigma}.
 \label{65}
 \eeq
\noindent Here the parameters $\eta^\sigma$ corresponds to
coordinates of the fiber (for instance, we can take the value
of the function $x(t)$ at $t=t_0$). It is natural to consider
$\eta^\sigma$  as coordinates of the standard fiber of the bundle
$\xi$. The  variational equations on a solution
$\hat x(\cdot), \hat q^*(\cdot)$ have the form

\beq
 \left \{
\begin{array}{l}
\frac {\partial V}{\partial t^m} = - \hat H_{q^*x}V -
 \hat H_{q^*q^*}U \\

 \frac {\partial U^m}{\partial t^m} = \hat H_{xx}V +
 \hat H_{xq^*}U \\
\end{array}
\right .
 \label{66}
 \eeq

\begin{assumption}
Let $U,V$ be a solution of variational equations~(\ref{66}) for
the canonical system of the Euler equations, defined on a domain
$\goth N$ of the space $t$.

Suppose that the matrix $V$ is defined and invertible in $\goth N$,
\end{assumption}

Consider the matrices $W^m, \; m=1...n$, which are defined by the
formulas $W^m_{\mu \rho} =  U^m_{\mu \sigma} (V^\rho
_\sigma )^{-1}$ on the domain $\goth N$. By differentiating $W^m$
we get

$$
\frac {\partial }{\partial t^m} (UV^{-1}) = \hat H_{xx}VV^{-1} +
\hat H_{xq^*}UV^{-1} + UV^{-1}(\hat H_{q^*x}V + \hat
H_{q^*q^*}U)V^{-1}.
$$

By substituting here the definition of $W$ we get the Riccati
equation in partial derivatives~(\ref{1}) for the second
variation of the functional~(\ref{1}) that corresponds to the
Hamiltonian $H$:

\beq
 \frac {\partial W^m}{\partial t^m} = \hat H_{xx} + \hat H_{xq^*}W +
 W\hat H_{q^*x} + W \hat H_{q^*q^*}W .
 \label{67}
\eeq

This equation defines the changing of $k$-Lagrangian planes
with the coordinates $W^m_{\rho \sigma}(t)$ along the extremal surface
 $\hat x(t)$ (See.~\cite{Ze}, \cite{Zel}).

\smallskip

\begin{theorem}
Let $\cal L$ be a smooth foliation on a fibre space $\xi$
generated by an invariant integral $\goth S$, and its fibers
have diffeomorphic projection on a domain $\goth N$ of the space
of variables $t$. Let the assumption 1 be fulfilled.

Then the functions
\beq
W^m_{\rho \sigma} = \frac {\partial }{\partial
x^\rho}
 \left \{\sum_{K\ni (m,l)} \adj^l_m
 \left \| \frac {\partial S^i}{\partial t^j} +
 \frac {\partial S^i }{\partial x^\alpha}
\frac {\partial x^\alpha}{\partial t^j} \right \|_K^t
 \right \} \frac {\partial S^l}{\partial x^\mu}
 \label{68}
\eeq

defines the solution to the Riccati equation~(\ref{67}).
\end{theorem}

\Proof.

The functions
$$
U^m_{\rho \beta} = \frac {n\partial (q^*)^m_\rho}{k\partial
\eta^\beta}
$$
define derivatives from $q^*$ with respect to coordinates of the
standard fiber. The matrix $V$ gives the operator of
differentiation of the mapping from the moving fiber into the
standard one of the space $\xi$. In view of assumption 1, this
mapping is invertible and its invers is defined by the matrix

$$
(V^{-1})^\beta _\sigma = \frac {\partial \eta^\beta}{\partial
x^\sigma}.
$$

The composition $UV^{-1}$ gives

$$
\frac {\partial [(q^*)^t]^m_\rho }{\partial \eta^\beta } \frac
{\partial \eta^\beta}{\partial x^\sigma} = \frac {\partial
[(q^*)^t]^m_\rho } {\partial x^\sigma}.
$$

In view of the formula~(\ref{31}) we get
$$
[(q^*)^t]^i_\alpha (t,x,q) =
 \frac {\partial \tilde {\Delta}}{\partial q_i^\alpha}
$$

Functions $\tilde {\Delta}$ и $\tilde {\goth S}$ are the
different form of the same integrand. Hence
$$
 \frac {\partial \tilde {\Delta}}{\partial q_i^\alpha} =
 \frac {\partial \tilde {\goth S}}{\partial q_i^\alpha}.
$$

The theorem is proved.

$\Box$

\smallskip

\begin{corollary}
In the case $(k=1)$ i.e. for the Weyl transform~\cite{W})
we obtain from~(\ref{68}) that the Hessian of the action-vector

$$
 W = UV^{-1} = \frac {\partial^2 S^m}{\partial x^\rho \partial
x^\sigma}.
$$
is the solution to the Riccati equation~(\ref{67})
\end{corollary}

Let us show that~(\ref{68})  defines
the connection on the fibre space $\xi$. With this in mind,
it is convenient to return from the canonical variables to
the Lagrangian ones. Since the transform  ${\cal Z}$ is
involutori, we have again to make the transform  ${\cal Z}$
of the quadratic approximation of the Hamiltonian that generates
the variational equation~(\ref{66})

\beq
 \frac {1}{2} ((H_{q^*q^*}U,U) + 2(H_{q^*x}U,V) +  (H_{xx}V,V)).
 \label{69}
\eeq

As a new variable we take the derivative of~(\ref{69}) with respect to $V$.

$$
H_{q^*q^*}U + H_{q^*x}V.
$$

Then $W$ transforms into

 \beq
\goth Y^{\alpha }_{i\beta } = (H_{q^*q^*}U + H_{q^*x}V)V^{-1} =
H_{q^*q^*}W + H_{q^*x}.
 \label{70}
 \eeq

As a differential form of connection we consider

\beq
 \zeta^\alpha  = dx^\alpha - \goth Y^\alpha _{i\beta}x^\beta dt^i.
 \label{71}
\eeq

The operator of covariant differentiation associated with the
form~(\ref{71}) is

\beq
 \nabla_{v^i}x^\alpha  = v^i\frac {\partial x^\alpha}{\partial t^i} -
 v^i\goth Y^\alpha _{i\beta}x^\beta .
 \label{72}
\eeq

The operator of projection of tangent vectors $(dt,dx)$ at a point
$(t,x)$ of the bundle $\xi$ on the fiber has the form

$$
(dt^i,dx^\alpha ) \mapsto (0,dx^\alpha -  \goth Y^\alpha
_{i\beta}x^\beta dt^i).
$$

Horizontal vectors of the connection $\nabla$ are vectors
$(dt,dx)$ that belongs to the kernel of that operator:
$dx - \goth Y_ixdt^i = 0$. Hence, the horizontal component
of the vector $(dt,dx)$ is

\beq
 (dt,\goth Y_ixdt^i).
\label{73}
 \eeq

A commutator of matrices $A$ and $B$ we will denote, as is customary,
by $[A,B]$.

\smallskip

\begin{theorem}

The tensor of curvature of the connection $\nabla$ equals

\beq
 \goth R  = \frac {\partial \goth Y_i}{\partial t^j} -
  \frac {\partial \goth Y_j}{\partial t^i} -
 [\goth Y_i \goth Y_j].
 \label{74}
\eeq
\end{theorem}

\Proof.

Let us recall that the exterior covariant derivative of
the form of connection $\zeta$ is called the form of
curvature of the given connection~\cite{K}. The exterior
covariant derivative $D\zeta$ is the value of the exterior
derivative  $d\zeta$ on the horizontal components of vectors
$(dt,dx)$. Let us calculate it

$$
d\zeta = - \goth Y_\alpha ^{m \beta}dx^\beta \wedge dt^m + \left (
\frac {\partial}{\partial t^j}\goth Y^\alpha_{i\beta}dt^i \wedge
dt^j - \frac {\partial}{\partial t^i}\goth Y^\alpha_{j\beta}dt^i
\wedge dt^j
 \right ) .
$$

By substitution of the horizontal component of the vector
$(dt,dx)$ we obtain~(\ref{74}).

$\Box$

\smallskip

A connection with the zero curvature is called flat.

\begin{theorem}

Let the assumption 1 be fulfilled. Then the connection $\nabla$
that was generated by the field of extremals is flat.

\end{theorem}

\Proof.

Let us rearrange the formula~(\ref{70}).

\beq
 \goth Y  = (H_{q^*q^*}U + H_{q^*x}V)V^{-1} = -
  \frac {\partial V}{\partial t^m}V^{-1}.
 \label{75}
\eeq

Substitute it in~(\ref{74}). We get

$$
\begin{array}{c}
- \frac {\partial^2 V}{\partial t^i\partial t^j}V^{-1} + \frac
{\partial V}{\partial t^i}V^{-1}\frac {\partial V}{\partial
t^j}V^{-1} - \left ( - \frac {\partial^2 V}{\partial t^j\partial
t^i}V^{-1} + \frac {\partial V}{\partial t^j}V^{-1}\frac {\partial
V}{\partial t^i}V^{-1} \right ) -  \\

 \frac {\partial V}{\partial t^i}V^{-1}\frac {\partial V}{\partial
t^j}V^{-1} +  \frac {\partial V}{\partial t^j}V^{-1}\frac
{\partial V}{\partial t^i}V^{-1} = 0.
\end{array}
$$

$\Box$

\smallskip

Hence, it was shown that fields of extremals define in the
given chart the flat curvature. It is because that the
existence of a field leads to horizontal integrable
distribution of planes. Its integral surfaces give the
full system of horizontal sections of the bundle $\xi$.

\bigskip

\section{Examples}

\bigskip

As an example let us consider the standard Hopf bundle $S^3\to S^2$.
Here $S^3$ is realized in the space $\C^2$ as a unit sphere:
$(z_1,z_2;\, w_1,w_2),\; |z|^2+|w|^2=1$. Fibers of the Hopf bundle
are defined as big circles $\{e^{i\phi}z, e^{i\phi}w\}$ passing through
each point $(z,w)$. The central projection $\pi$ from the center of
the sphere onto the tangent plane $P$ at the point $(1,0; \, 0,0)$
(at the north pole of $S^3$) turns $P$ into a 3-dimensional projective space
$\R\P^3$ with the coordinates that corresponds to the three last
coordinates of the plane $P$. It may be considered as a chart
$\goth A$ on the north hemisphere of $\S^3$
The metric on the space $\R\P^3$.
induced by the projection $\pi$, is

$$
ds^2 = (1+\eta^2 + \zeta_1^2 + \zeta_2^2)^{-2}
 \left (
(d\eta)^2 + (d\zeta_1)^2) +(d\zeta_2)^2.
 \right ).
$$

Fibers of the Hopf bundle projects on straight lines (one of two
families of rectilinear generators of the set of hyperboloid of
one sheet --- projections of torus on $S^3$). So, we have on $P$
the family of geodesic $\goth P$, and through each point passes
one and only one geodesic. Each geodesic gives the absolute
minimum of the length among all curves lying in the chart $\goth A$.
Let us show that there is no action-function that
simultaneously synchronizes all the extremals of the set $\goth P$.
Indeed, fibers of the Hopf bundle are obtained by the simultaneous
rotation about the same angle $\varphi$ in planes $z$ and $w$. The
point $(1,\eta; \, \zeta_1, \zeta_2)$ passes into the point

$$
A=((\cos \varphi + \eta \sin \varphi),(\sin \varphi - \eta \cos
\varphi),(\zeta_1\cos \varphi + \zeta_2\sin \varphi),(\zeta_1\sin
\varphi - \zeta_2\cos \varphi)).
$$

The image of the $\pi$ projection of the point $A$ is obtained by
normalization --- dividing by the first coordinate. So, we get

$$
 \left (
\frac {\sin \varphi - \eta \cos \varphi}{\cos \varphi + \eta \sin
\varphi}, \frac {\zeta_1\cos \varphi + \zeta_2\sin \varphi}{\cos
\varphi + \eta \sin \varphi}, \frac {\zeta_1\sin \varphi -
\zeta_2\cos \varphi}{\cos \varphi + \eta \sin \varphi}
 \right ).
$$

Tangent vectors $X$ to fibers of the Hopf bundle are obtained
by differentiation with respect to $\varphi$ and putting $\varphi =0$.
Hence, $X = ((1+ \eta^2), (-\eta \zeta_1 +\zeta_2), (\zeta_1
+ \eta \zeta_2))$. The transversality condition for $n=1$ is the
orthogonality condition in the metric induced on $P$ by the
projection $\pi$. Orthogonal planes to vectors $X$
correspond  to the space of zeroes of the differential form
$\xi = (1+\eta^2)d\eta + (-\eta \zeta_1+\zeta_2)d\zeta_1 +
(\zeta_1+\eta \zeta_2)d\zeta_2$. However, $\xi \wedge d\xi \ne 0$
and the form $\xi$ is not integrable. It is the obstacle to
design action-function. Consequently, $\goth P$ does not generate
a field of extremals. The corresponding manifold is not a
Lagrangian one.

\bigskip

The similar situation takes place for the quaternary Hopf
bundle $S^7 \to S^4$ \cite{N}. It is realized as a unit sphere
 $\{(z,w), \; \|z\|^2 + \|w\|^2 = 1\}$ in the 2-dimensional
quaternary space $\H^2$. A fiber passing through a point $(z,w)$
is defined as a set of points $\{\sigma z, \sigma w\}$ where
$\sigma$ runs over the unit quaternary sphere $\|\sigma\|=1$. It is
easy to verify that each fiber is the central section of the
sphere $S^7$. Through each point of the sphere passes one and
only one of these fibers. The projection $\pi$ from the center
of the sphere on the tangent plane at the north pole
$(1,0,0,0; 0,0,0,0)$ of the sphere turns $P$ into 7-dimensional
projective space $\R\P^7$ with coordinates $(\eta_1,\eta_2,\eta_3;
 \zeta_0,\zeta_1,\zeta_2,\zeta_3)$ that corresponds to the seven
last coordinates of the plane $P$. The metric on the space $\R\P^7$
induced by projection $\pi$ is

$$
ds^2 = (1+(\eta)^2+(\zeta)^2)^{-2}((d\eta)^2+(d\zeta)^2).
$$

The projection $\pi$ transfers fibers into the set $\goth P$ of
3-dimensional planes. Through each point passes one and only one
plane. Each such plane gives the absolute minimum to the
functional of 3-dimensional volume in $\R\P^7$ in the class of
variations lying in the considered chart. The normalization
relative to the first coordinate defines the parametric equation
of the planes $\goth P$:

$$
\eta_i(\sigma ) = (\sigma z)_0^{-1}(\sigma z)_i, \; (i=1,2,3);
\quad \zeta_j(\sigma ) = (\sigma z)_0^{-1}(\sigma w)_j, \;
(j=0,1,2,3).
$$

Note that an element of $k$-dimensional volume is determined
by the length of $k$-dimensional multivector in metric being
the tensor $k$-power of the metric of ambient space~\cite{R}.
The integrand in the situation in question is

$$
f = (1 + \eta^2 + \zeta^2)^{-3}\sqrt{\sum_{|I|+|J|=3}
 \left | \frac {D(\eta^I\zeta^J)}{D(\sigma)}
 \right |^2}.
$$

The function $f$ is not a convex one and the Weyl construction
is inapplicable. However, the function $f$ is convex as a
function of Jacobians

$$
 \left | \frac {D(\eta^I\zeta^J)}{D(\sigma)}
 \right |.
$$
Hence, we can apply the developed theory taking $k=3$.
Let us use the transversality condition from the section 5. It
can be verified that the distribution of planes that are transverse
to fibers of the quaternary Hopf bundle is not an integrable one.
The family $\goth P$ does  not give the field of extremals and
the corresponding manifold is not a 3-Lagrangian one.

Calculations in this case are much more cumbersome and we do
not give it here.

It is possible to calculate the differential form of connection
for the distribution of normal planes to fibers of the
Hopf bundle and the corresponding form of curvature. It may be presumed
that we obtain the Chern class of the Hopf bundle.

It would be interesting to consider exotic 7-dimensional
Milnor's spheres for which fibers are 3-dimensional spheres
that are obtained by the action of unit quaternary $\sigma$
on $S^7 \subset \H^2$ using the formula $\{\sigma z, \sigma^h
w\sigma^j\}$, where $h+j=1$ (see. \cite{Mi}). Here is an
additional difficulty connected with the finding of the
explicit expression for metric on these spheres.


\begin{thebibliography}{99}

\bibitem{N} B.A.Dubrovin, S.P. Novikov, A.T.Fomenko.
"The modern Geometry". Moscow. Nauka. (1979)(in Russian)
\bibitem{Ze} M.I.Zelikin. "Homogenious spaces and the Riccati
equation in the calculus of variations." Moscow. Faktorial.
(1998)(in Russian)
\bibitem{Zel} M.I.Zelikin. "Hessian of solutions to Hamilton-Jacobi
equation in the theory of extremal problems" Matematicheskii
Sbornik. 195:6 (2004) pp.819-831
\bibitem{I}A.D.Ioffe, V.M.Tikhomirov, "The theory of extremal problems."
Moscow. Nauka.(1974)(in Russian)
\bibitem{K} Sh.Kobajashi,K.Nomizu. "Foundations of differential geometry"
Interscience publishers. New York. London.(1863)
\bibitem{M} J.Milnor,J.Stasheff. "Characteristic classes".
Princeton. New Jersy.(1974)
\bibitem{R} P.K.Rashevskii. "Riemannian geometry and tensor analise".
Moscow. Nauka.(1967). (in Russian)
\bibitem{Fo} A.T.Fomenko. "Variational methods in topology."
Moscow. Nauka. (1982). (in Russian)
\bibitem{Ba} J.M.Ball, F.Murat.  "$W^{1,p}$-quasiconvexity and
variational problems for multiple integrals". J.Func.Anal. V.58
(1984), pp.225-255
\bibitem{Ba1} J.M.Ball. "Remarks on the paper 'Basic calculus of
variations'", Pacific J. of Math.. Vol.116, № 1, (1985), pp.7-9
\bibitem{B} H.B\"orner. \"Uber die Legendrische Bedingung und die
Feldtheorien in der Variationsrechnung der mehrfachen Integrale.
Math.Zeitschr. 46. (1940)
\bibitem{C} C.Caratheodory. \"Uber die Variationsrechnung bei
mehrfachen Integralen". Acta Szeged, 4 (1929) S.193-216
\bibitem{Cl} A.Clebsch.  \"Uber die zweite Variation vielfacher
Integrale. J.Reine Angew.Mathematik", Bd 56, (1859) S.122-148
\bibitem{Du} J.J.Duistermaat. "On the Morse Index in Variational
Calculus".  Advances in Math. 21 (1976) pp.173-195
\bibitem{D} Th. De Donde. "Theorie invariantive du Calcul des
Variations". Paris, Gauthier-Villars (1935)
\bibitem{Den} Dennemeyer. "Conjugate Surfaces for Multiple Integral
Problems." Pacific J. Math. V.30 № 3,(1969), pp.621-638
\bibitem{F} P.Funk. "Variationsrechnung und ihre Anwendungen in
Physik und Technik". Berlin-G\"ottingen-Heidelberg.
Springer-Verlag. (1962)
\bibitem{G} M.Giaquinta, S.Hildebrandt. "Calculus of Variarions."
Springer-Verlag, Berlin-Heidelberg-New York,(1996)
\bibitem{Gi} E.Giusti. "Direct Methods in the Calculus of
Variations". World Scientific. New-Jersey-London- Singapore- Hong
Kong (2005)
\bibitem{H} J.Hadamard. "Sur quelque question de calcul des
variations". Bull.Soc.Math.France. V.33. (1905).p.77-80
\bibitem{L} J.T.Lepage. "Sur le champs geodesiquedes
integrales multiples". Bull.Acad.Roy.Belg. Class de sciences.
27.(1936),p.27-46
\bibitem{Mi} J.Milnor.  "On manifolds homeomorphic to the
7-sphere". Ann.Math. 64 (1956), p.399-405
\bibitem{Mo} C.B.Morrey. "Multiple Integrals in the Calculus of
Variations". Springer Verlag, Vol.130, (1966)
\bibitem{S} D.Serre. "Condition de Legendre-Hadamard; espaces de
matrices des rang $\ne 1$".  C.R.Acad.Sci.Paris Ser.I Math., 293
(1981),p.23-26
\bibitem{SE} D.Serre. "Formes quadratiques et calcul des
variations". J.Math.Pures et Appl. 62 (1983), p.177-196
\bibitem{Sil} E.Silverman. "A sufficient condition for the lower
semicontinuity of parametric integrals". Trans.Amer.Math.Soc. 167,
(1972), pp.465-469
\bibitem{Si} J.Simons. "Minimal varieties in Rimannian manifolds".
Ann.of Math. 88, (1969), pp.62-105
\bibitem{Sm} S.Smale. "On the Morse index theorem". J.Math.Mech.
14 (1965), pp.1049-1056
\bibitem{Sw} R.C.Swanson. "Fredholm intersection theory and elliptic
boundary deformation problems". J.Diff.Eq. V28, № 2.(1978),
pp.189-219
\bibitem{T} F.J.Terpstra. "Die Darstellung biquadratischer Formen
als Summen von Quadraten mit Anwendungen auf die
Varitionsrechnung". Math Ann.116 (1938) S.166-180
\bibitem{U} K.\"Uhlenbeck. "The Morse index theorem in Hilbert Space"
 J.Diff.Geometry 8 (1973), pp.555-564
\bibitem{V} L.Van Hove. "Sur l'extension de la condition de
Legendre du calcul des Variations aux integrales multiples a
plusiers fonctions inconnues". Nederl.Akad.Wetensch. 50 (1947)
p.18-23
\bibitem{W} H.Weyl. "Geodesic fields in the calculus of variations
for multiple integrals". Ann.Math. 36. (1935) p.607-629
\bibitem{Z} M.I.Zelikin. "Control Theory and optimization. I".
Vol.86. Springer-Verlag, Berlin-Heidelberg-New York,(2000)

\end{thebibliography}
\end{document}